\documentclass[12pt]{article}
\usepackage{amsmath,amssymb,amsfonts,verbatim}
\input{amssym.def}\input{amssym}

\newcounter{theorem}\makeatletter

\newtheorem{theorem}{Theorem}\newtheorem{prop}[theorem]{Proposition}\newtheorem{lem}[theorem]{Lemma}
\newtheorem{def-lem}[theorem]{Definition-Lemma}\newtheorem{def-prop}[theorem]{Definition-Proposition}
\newtheorem{cor}[theorem]{Corollary}

\def\a{\alpha}\def\Ar{{\rm A}}\def\Ab{\bar{\rm A}}\def\ad{\operatorname{ad}}\let\al=\a\def\bb{\mathfrak{b}}
\def\ba{\begin{eqnarray}}\def\ea{\end{eqnarray}}\def\be{\begin{equation}}\def\ee{\end{equation}}\def\Cprime{C'}\def\com{\ts,\hskip-.5pt}
\def\ds{\displaystyle}\def\De{\Delta}\def\de{\delta}\def\f{\mathfrak{k}}
\def\h{\mathfrak{h}}\def\hs{\mathring{z}}
\def\g{\mathfrak{g}}\def\sl{\mathfrak{sl}}\def\gl{\mathfrak{gl}}\def\I{\mathrm{I}}\def\Ib{\tts\overline{\nns\rm I}_+\tts}
\def\Jb{\,\overline{\nns\rm I}_-\ts}\def\k{\mathfrak{k}} \def\kar{a}\def\la{\lambda}\def\lcd{\ts,\ldots,}
\def\mult{\diamond}\let\mmult=\mult\def\n{\mathfrak{n}}\def\nminus{\mathfrak{n}_-}
\def\nplus{\mathfrak{n}_+}\def\nns{\hskip-.5pt}\def\p{\mathfrak{p}}
\def\pp{\widetilde{\phantom{\,}p\phantom{\,}\phantom{\ }}\phantom{\!}\phantom{\!}\phantom{\!}}\def\q{\check{\mathrm q}}\def\Q{\mathrm{Q}}
\let\rf=\r\def\s{z}\def\S{\mathrm {S}}\def\si{\sigma}\def\sy{\acute{\sigma}}
\def\Sym{\operatorname{W}}\def\th{\mathring{h}}\def\tphi{A}\def\tpsi{A'}\def\ttphi{B}\def\ttpsi{B'}
\def\ts{\hskip1pt}\def\tts{\hskip.5pt}\def\tt{\mathring{t}}\def\U{\mathrm{U}}\def\Uh{\overline{\U}(\h)}\def\UUh{{\cal{D}}(\h)}
\def\V{\mathrm {V}}\def\ve{\varepsilon}\def\Z{\mathrm{Z}}\def\z{z}\def\Zgk{\Z(\g,\k)}
\newcommand{\ZZ}{{\mathbb Z}}

\newcommand\qedsymbol{\hbox{\rlap{$\sqcap$}$\sqcup$}}\newcommand\qed{\relax\ifmmode\else\unskip\quad\fi\qedsymbol}

\begin{document}

\begin{center}
{\Large\bf Diagonal reduction algebras of $\gl$ type}

\bigskip
{\bf  S. Khoroshkin$^{\circ}$ \ \ and \, O. Ogievetsky$^{\star}$\footnote{On leave of absence
from P.N. Lebedev Physical Institute, Theoretical Department, Leninsky prospekt 53, 119991 Moscow, Russia}
}\medskip\\
$^\circ${\it Institute of Theoretical and Experimental Physics, 117259 Moscow, Russia}\smallskip\\
$^\star${\it Centre de Physique Th\'eorique\footnote{Unit\'e Mixte de Recherche (UMR 6207) du CNRS et des Universit\'es Aix--Marseille I,
Aix--Marseille II et du Sud Toulon -- Var; laboratoire affili\'e \`a la FRUMAM (FR 2291)}, Luminy, 13288 Marseille, France}
\end{center}

\medskip
\setcounter{equation}{0}
\begin{abstract}\noindent
Several general properties, concerning reduction algebras -- rings of definition and algorithmic efficiency of
the set of ordering relations -- are discussed. For the reduction algebras, related to the diagonal embedding
of the Lie algebra $\gl_n$ into $\gl_n\oplus\gl_n$, we establish
a stabilization phenomenon and list the complete sets of defining relations.
\end{abstract}

\section{Introduction} Reduction algebras were introduced \cite{AST2,M} for a study of representations of a
Lie algebra with the help of the restriction to a subalgebra.

Let $\g$ be a Lie algebra, $\f\subset \g$ its reductive Lie subalgebra; that is,
the adjoint action of $\k$ on $\g$ is completely reducible (in particular,
$\k$ is reductive).
Suppose $\k$ is given with  a triangular decomposition
\be\label{intro1} \k=\nminus+\h+\nplus\ .\ee
Denote by $\I_+$ the left ideal of $\Ar :=\U(\g)$ generated by elements of $\nplus$, $\I_+ :=\Ar\nplus$. Then the reduction algebra $\S(\g,\k)$, related to the pair
$(\g,\k)$, is defined as the quotient $\mathrm{Norm}(\I_+)/\I_+$ of the normalizer  of the  ideal $\I_+$ over $\I_+$ (one should keep in mind that the notation $\S(\g,\k)$
is abbreviated: the data needed for the definiton of the reduction algebra includes, in addition to the pair $(\g,\k)$, the triangular decomposition \rf{intro1}).
The space $\S(\g,\k)$ is equipped with a natural structure of the associative algebra.  By definition, for any $\g$-module $V$ the space $V^{\nplus}$ of vectors,
annihilated by $\nplus$, is a module over $\S(\g,\k)$. If $V$ decomposes, as an $\f$-module, into a direct sum of irreducible $\k$-modules $V_i$ with finite-
dimensional multiplicities, then the $\g$-module structure on $V$ can be uniquely restored from the $\S(\g,\k)$-module structure on $V^{\nplus}$.

The reduction algebra simplifies  after the localization over the multiplicative set generated by elements $h_\gamma+k$, where $\gamma$ ranges through the
set of roots of $\k$, $k\in\ZZ$; here $h_\gamma$ is the coroot corresponding to $\gamma$. Let $\Uh$ be the localization of the universal enveloping algebra
$\U(\h)$ of the Cartan sub-algebra $\h$ of $\k$ over the above multiplicative set. The localized reduction algebra $\Z(\g,\k)$ is an algebra over the commutative
ring $\Uh$; the principal part of the defining relations is quadratic but the relations may contain linear terms or degree 0 terms, see \cite{Zh,KO}. Besides, the reduction
algebra admits another description as a (localized) double coset space
$\Ar/(\Ar \nplus +\nminus\Ar)$
endowed with the multiplication map defined with the help of the extremal projector \cite{KO} of Asherova--Smirnov--Tolstoy \cite{AST}.

The general theory of reduction algebras \cite{Zh} provides the set of quadratic-linear-constant
ordering relations over $\UUh$, the field of fractions of  $\U (\h )$, for natural generators of any reduction algebra $\Z(\g,\f)$.  However there are two subtle questions concerning
these relations. The first question is:  are these ordering relations defined over the smaller ring $\Uh$?
Secondly, is it possible to use these ordering relations for an algorithmic ordering of polynomial expressions in the reduction algebra?
In the first part of the paper we give affirmative answers to these questions for any reduction algebra $\Zgk$.

The main theme of the second part of the paper is the special restriction problem, when $\g$ is the direct sum of two copies of the Lie algebra ${\gl_n}$
and $\f$ is the diagonally embedded ${\gl_n}$. The resulting reduction algebra we call {\em diagonal reduction} algebra of $\gl_n$ and denote by $\Z_n$.
A finite-dimensional irreducible module over $\g=\gl_n\oplus\gl_n$ is the tensor product of two irreducible $\gl_n$-modules; restricting the $\g$-module
to $\k$ we obtain the decomposition of the tensor product into the direct sum of irreducible $\gl_n$-modules. One of the main results of the paper is the explicit
description of the diagonal reduction algebra $\Z_n$. Some examples and applications of the diagonal reduction algebras are given in \cite{KO3}.

We present a list of defining relations for natural generators of $\Z_n$. The derivation of these relations uses heavily the Zhelobenko automorphisms \cite{KO}
of reduction algebras and is given in the work \cite{KO2}. In the present paper we formulate and prove the stabilization property of the algebras $\Z_n$.
The stabilization phenomenon provides a natural way of extending relations for $\Z_{n}$ to relations for $\Z_{n+1}$ ($\Z_{n}$ is not a subalgebra of $\Z_{n+1}$).
The stabilization principle is the second essential ingredient for the derivation of the set of defining relations.

We also prove that our list of defining relations is equivalent over $\Uh$ to the list of the canonical ordering relations. The proof is not difficult once we treat the
algebras over $\UUh$: the arguments for the equivalence are based on certain asymptotic considerations. The proof of the equivalence over $\Uh$ is more
delicate, it uses the stabilization phenomenon and calculations of certain determinants of Cauchy type.

\underline{\underline{}}\section{Reduction algebras related to a reductive pair}
\label{section2new}

Let $\g$ be a finite-dimensional Lie algebra and $\k\subset\g$ its reductive subalgebra. Assume that
the embedding $\k\subset\g$ is also reductive, that is the adjoint action of $\k$ in $\g$ is semi-simple. Let
$\p$ be an $\ad_\k$-invariant complement of $\k$ in $\g$. Choose a triangular decomposition \rf{intro1} of Lie algebra $\k$; here
 $\h$ is a Cartan subalgebra of $\k$ while $\nplus$ and $\nminus$ are nilradicals of two
opposite Borel subalgebras $\bb_\pm\subset\k$.
Let $\De\in\h^*$ be
the root system of $\f\ts$. The subsets of $\De$ consisting of the positive
and negative roots will be denoted by $\De_+$ and $\De_-$
respectively.
Let $\Q$ be the root lattice,
$\Q:=\{\gamma\in\h^*\,|\,\gamma=\sum_{\alpha\in\Delta_+, n_\a\in\ZZ}n_\a \a\}$. It contains the
positive cone $\Q_+$,
\be\Q_+:=\{\gamma\in\h^*\,|\,\gamma=\!\!\!\sum_{\alpha\in\Delta_+, n_\a\in\ZZ,
n_\a\geq 0}n_\a \a\}\ .\label{cqp}\ee
{}For $\lambda,\mu\in\h^*$, the notation
\be\lambda>\mu\label{paor}\ee
means that the difference $\lambda-\mu$ belongs to $\Q_+$, $\lambda-\mu\in\Q_+$.
This is a partial order in $\h^*$.

Let $\Sym$ be the Weyl group of the root system $\De\ts$. Let
$\si_1\lcd\si_r\in\Sym$ be the reflections in $\h^*$ corresponding
to the simple roots $\al_1\lcd\al_r\ts$. We also use  the
induced action of the Weyl group $\Sym$ on the vector space $\h\ts$.
It is defined by setting $\la\ts(\si(H))=\si^{-1}(\la)(H)$ for all
$\si\in\Sym$, $H\in\h$ and $\la\in\h^*$. We assume that this action
is extended to the action of a cover of the group $\Sym$ by automorphisms of
the Lie algebra $\g$. In other words, there are automorphisms $\sy_i:\g\to\g$
which satisfy the same braid group relations as $\si_i$, preserve the subspaces $\h$ and $\k$,
and coincide with $\si_i$ being restricted to $\h$. We denote by the same symbols
the canonical extensions of $\sy_i$ to automorphisms of $\U(\g)$.

Let $\rho$ be the
half-sum of the positive roots of $\k$. Then the \textit{shifted action\/}
$\circ$ of the group $\Sym$ on the vector space $\h^*$ is defined by
setting
\be\label{shifted}
\si\circ\la=\si(\la+\rho)-\rho\ts.
\ee
With the help of \rf{shifted} we induce the action $\circ$ of $\Sym$ on the commutative
algebra $\U(\h)\ts$
by regarding the elements of this algebra as
polynomial functions on $\h^\ast$. In particular, then
$(\si\circ H)(\la)=H(\si^{-1}\circ\lambda)$ for $H\in\h\ts$.

{}For each $i=1\lcd r$ let $h_{\a_i}=\al_i^\vee\in\h$ be the coroot vector
corresponding to the simple root $\al_i\ts$, so that the value
$\al_j\ts(H_i)$ equals the $(i\com j)$ entry $\kar_{ij}$ of the
Cartan matrix $\kar\ts$ of $\k$. Here $h_{\a_i}$ belongs to the semi-simple
part of $\f\ts$. Let $e_{\a_i}\in\nplus$ and $e_{-\a_i}\in\nminus$ be the
Chevalley generators of that subalgebra corresponding to the roots
$\al_i$ and $-\al_i\ts$ so that
$$
[e_{\a_i}\com e_{-\a_j}]=\de_{ij}\,h_{\al_i}\,, \quad [h_{\al_i},e_{\al_j}]=
\kar_{ij}\,e_{\al_j}\,,
\quad [h_{\al_i},e_{-\al_j}]=-\kar_{ij}\,e_{-\al_j}\, .
$$
{}For each $\al\in\De$ let $h_\al=\al^\vee\in\h$ be the corresponding
coroot vector. Denote by $\Uh$ the ring of fractions of the
commutative algebra $\U(\h)$ relative to the set of denominators
\begin{equation}
\label{M2} \{\,h_\al+l\ |\ \al\in\Delta\ts,\ l\in\ZZ\,\ts\}\,.
\end{equation}
The elements of this ring can also be regarded as rational functions
on the vector space $\h^\ast\ts$. The elements of
$\U(\h)\subset\,\overline{\!\U(\h)\!\!\!}\,\,\,$ are then regarded
as polynomial functions on $\h^\ast\ts$. Let
$\overline{\!\U(\k)\!\!\!}\,\,\,\subset \Ab= \overline{\!\U(\g)\!\!\!}\,\,\,$ be the rings of
fractions of the algebras $\U(\k)$ and $\Ar=\U(\g)$ relative to the set of denominators
\eqref{M2}. These rings are well defined, because both $\U(\f)$ and $\U(\g)$ satisfy the Ore
condition relative to \eqref{M2}. Since $\si_i$ preserve the set of denominators \rf{M2},
the automorphisms $\sy_i$ admit a natural extension to $\Ab$.

Define $\Z(\g,\k)$ to be the double coset space of $\Ab$ by its left ideal $\Ib:=\Ab\nplus$,
 generated by elements of $\nplus$, and the right ideal $\Jb:=\nminus\Ab$,
generated by elements of $\nminus$, $\Z(\g,\k):=\Ab/(\Ib+\Jb)$.
The space $\Z(\g,\k)$ is an associative algebra with respect to the multiplication map
\begin{equation}\label{not5a}a\mult b:=aP b\ .\end{equation}
Here $P$ is the extremal projector \cite{AST} of the Lie algebra $\k$ corresponding to the triangular
decomposition \rf{intro1}. We call $\Z(\g,\k)$
the {\textit{reduction}} algebra associated to the pair $(\g,\k)$.
The assignment  $\,x\mapsto x\mod \Ib+\Jb\,$ establishes an injective homomorphism of the algebra
$\S(\g,\k)$ (see Introduction for the definition)
to $\Zgk$, see \cite{KO}. Moreover, the
localization of the image of $\S(\g,\k)$
with respect to $\Uh$ coincides with $\Zgk$.

The algebra $\Z(\g,\k)$ can be equipped with the action of
Zhelobenko automorphisms \cite{KO}.  Denote by $\q_{i}$  the Zhelobenko automorphism
$\q_i:\Zgk\to\Zgk$  corresponding to the simple root $\alpha_i$, $i=1\lcd r$.
It is defined as follows \cite{KO}. First we define a map $\q_i:\Ar\to \Ab/\Ib$ by
\begin{equation}\label{not7}\q_i(x):=\sum_{k\geq 0}\frac{(-1)^k}{k!}\hat{e}_{\a_i}^k(\sy_i(x))
e_{-\a_i}^k\
\ds\prod_{j=1}^k(h_{\a_i}-j+1)^{-1}
\quad \mod\Ib\ .\end{equation}
Here $\hat{x}$ stands for the adjoint action of the element $x$, so that $\hat{x}(y)=xy-yx$ for
$x\in\k$ and $y\in\Ab$.
The operator $\q_i$ has the property
\be\label{not2}\q_i(hx)=(\si_i\circ h)\q_i(x)\ee
for any $x\in\Ar$ and  $h\in\h$; $\si\circ h$ is defined in \rf{shifted}. With the help
of (\ref{not2}), the map $\q_i$ can be extended
to the map (denoted by the same symbol) $\q_i:\Ab\to \Ab/\Ib$ by the setting
$\q_i(\phi x)=(\si_i\circ \phi )\q_i(x)$ for any $x\in\Ar$ and
$\phi\in\Uh$. One can further prove that $\q_i(\Ib)=0$ and $\q_i(\Jb)\subset (\Jb+\Ib)/\Ib$,
so that $\q_i$ can be viewed as a linear operator $\q_i:\Zgk\to\Zgk$. Due to \cite{KO},
this is an algebra automorphism, satisfying \rf{not2}.
The operators $\q_i$ satisfy the same braid group relations as $\si_i$ and
the inversion relation \cite{KO}:
\begin{equation}\label{invr}\q_i^2(x)=(h_{\a_i}+1)^{-1}\ \sy_i^2(x)\
(h_{\a_i}+1)\ ,\qquad x\in\Zgk\ .\end{equation}

Let $\p$ be an $\ad_\k$-invariant complement of $\k$ in $\g$, as above.
Choose a linear basis $\{p_K\}$ of $\p$ and
equip it with a total order $\prec$. For an arbitrary element $a\in\Ab$ let $\widetilde{a}$
be its image in the reduction algebra;
in particular, $\pp_K$ is the image in $\Zgk$ of the basic vector $p_K\in\p$.

\vskip .2cm
The general theory of reduction algebras, see \cite{Zh} for the statements (a)-(c), says:

\begin{itemize}
\item[(a)] Since $\h$ normalizes both $\nplus$ and $\nminus$, the algebra $\Zgk$ is a
$\Uh$-bimodule with respect to the multiplication by elements of $\Uh$. It is
free as a left $\Uh$-module and as a right $\Uh$-module. As a generating
 (over $\Uh$) subspace one can take a projection of the space $\mathrm{S}(\p)$ of
symmetric tensors on $\p$ to $\Zgk$, that is a subspace of $\Zgk$, formed by linear combinations
of images of the powers $p^\nu$, where $p\in\p$ and $\nu\geq 0$.

\item[(b)] Assignments  $\deg (\widetilde{\hspace{.05cm}X\hspace{.05cm}})=l$ for the image
 of any product of $l$ elements from $\p$, $X=p_{K_1}p_{K_2}\cdots p_{K_l}$,
and  $\deg (Y)=0$ for any $Y\in\Uh$ define the structure of a
filtered algebra on $\Zgk$. The subspace $\Zgk^{(k)}$ of elements of
degree not greater than $k$ is a free left $\Uh$-module and a free
right $\Uh$-module, with a generating subspace formed by linear
combinations of images of the powers $p^\nu$, where $p\in\p$ and
$k\geq\nu\geq 0$.

\item[(c)] In the sequel we will choose for $\{p_K\}$ a weight ordered basis; that is, each $p_K$
has a certain weight $\mu_K$,
\be [h,p_K]=\mu_K(h)p_K\label{weir}\ee
for all $h\in\h$. The total order $\prec$ will be compatible with the partial order $<$ on $\h^*$,
see \rf{paor}, in the sense that
$\mu_K<\mu_L\ \Rightarrow p_K\prec p_L\ .$
Then the images $\pp_{\bar{L}}$ of the monomials ($\bar{L}$ is understood as the multiindex)
\be p_{\bar{L}}:=p_{L_1}^{n_1}p_{L_2}^{n_2}\cdots p_{L_m}^{n_m},\qquad p_{L_1}\prec p_{L_2}\prec\ldots\prec
 p_{L_m}\
,\qquad k=n_1+\dots+n_m\ ,\label{inimp}\ee
in $\Zgk^{(k)}$ are linearly independent over $\Uh$ and their projections to  the quotient
$\Zgk^{(k)}/ \Zgk^{(k-1)}$ form a basis of the left $\Uh$-module
$\Zgk^{(k)}/ \Zgk^{(k-1)}$.
The structure constants of the algebra $\Zgk$
 in the basis $\{\pp_{\bar{L}}\}$  belong to the ring $\Uh$.

\end{itemize}
Choosing the PBW basis of $\Ar$ induced by any ordered basis of $\k +\p$, which starts from a basis
in $\nminus$ and ends by a basis in $\nplus$, we see
that the statement about the monomials \rf{inimp} in (c) is valid without any condition
on the order $\prec$. However, the compatibility of the order $\prec$
with the partial order $<$ on $\h^*$ will be crucial for most of the statements below.

\begin{itemize}
\item[(d)] The algebra $\Zgk$ is the unital associative algebra, generated by $\Uh$ and all
$\{\pp_L\}$,
with the weight relations \rf{weir} and the ordering relations
\be\label{not3} \pp_I\mult \pp_J=\displaystyle{\sum_{K,L:p_K\preceq p_L}}{\mathrm{B}}_{IJKL}\,
\pp_K\mult\pp_L+\sum_M{\mathrm{C}}_{IJL}\pp_L+ {\mathrm{D}}_{IJ}\ ,\quad p_I\succ p_J\ , \ee
where ${\mathrm{B}}_{IJKL}$, ${\mathrm{C}}_{IJL}$  and ${\mathrm{D}}_{IJ}$ are certain
elements of $\Uh$.
\end{itemize}

Let  $\UUh$ be the field of fractions of the ring $\U (\h )$.
In \cite{Zh}, sections 4.2.3 - 4.2.4 and 6.1.5, it is proved that
the reduction algebra $\Zgk$ is generated by
the elements $\pp_L$ with the defining ordering relations \rf{not3} as an algebra over
 $\UUh$.
We shall now show that the statement (d) holds over the smaller ring $\Uh$; in other words, the
relations \rf{not3} are defined over $\Uh$ and the elements $\pp_L$ generate over $\Uh$
 the algebra $\Zgk$.

\vskip .2cm
We first prove that the structure constants ${\mathrm{B}}_{IJKL}$, ${\mathrm{C}}_{IJL}$
and ${\mathrm{D}}_{IJ}$ belong actually to $\Uh$. This fact can be understood with the help of the
factorized formula \cite{AST} for the extremal projector $P$. Indeed, decomposing the product,
we represent the projector $P$, after some reorderings, as a sum of terms
$\xi e_{-\gamma_1}\cdots e_{-\gamma_m}e_{\gamma'_{1}}\cdots e_{\gamma'_{m'}}$,
where $\xi\in\Uh$, $\gamma_1,\dots ,\gamma_m$ and
$\gamma'_1,\dots ,\gamma'_{m'}$ are positive roots of $\f$; the denominator of $\xi$ is a
product of linear factors of the form
$h_\gamma+\rho(h_\gamma)+\ell$,
where $\gamma$ is a positive root of $\f$ and $\ell$ a positive integer, $\ell >0$.
 We calculate the product $a\mult b$ in the following way.
 In the summand $a\xi e_{-\gamma_1}\cdots e_{-\gamma_m}e_{\gamma'_{1}}\cdots e_{\gamma'_{m'}}b$
 of $a\mult b$,
we move $\xi$ and all $e_{-\gamma}$'s to the left through $a$ by taking multiple commutators with $a$ and,
similarly, all $e_{\gamma'}$'s to the right through
$b$. Proceeding this way, we write
\be \pp_I\mult \pp_J={\mathrm{M}}_{IJKL}\widetilde{\phantom{\,}p_Kp_L}\label{matrM}\ee
(we recall that $\widetilde{a}$ denotes the image of an element $a\in\Ar$ in the reduction algebra)
where the (uniquely defined by the method of calculation)
matrix ${\mathrm{M}}$ with entries in
$\Uh$ has a triangular structure (even more is true: ${\mathrm{M}}_{IJKL}\neq 0$ $\Rightarrow$
$p_I\succ p_K$)
with 1's on the diagonal; denominators of
entries of the matrix ${\mathrm{M}}$ are of the form $h_\gamma
+\rho(h_\gamma)+ \pi(h_\gamma)+\ell$,
where $\pi$ is the weight, with respect to $\h$, of the
corresponding $p_I$ (the summand $\pi (h_\gamma )$ appeared when, in calculating
$\pp_I\mult \pp_J$
as above, we first moved
$\xi\in\Uh$ to the left through $\pp_I$; taking further multiple commutators, we do not change the
denominators any more). Take the formal (in the sense that for the moment we do not pay attention to possible dependencies between
$\widetilde{\phantom{\,}p_Ip_J}$ or between $\pp_K\mult \pp_L$ in the algebra) inverse: $\widetilde{\phantom{\,}p_Ip_J}=
{\mathrm{M}}_{IJKL}^{-1}\pp_K\mult \pp_L$;
the inverse matrix ${\mathrm{M}}^{-1}$ is
triangular as well, its entries are in $\Uh$ and it has 1's on the diagonal; the determinant of
${\mathrm{M}}$ is thus 1 and it follows that the above
described structure of denominators of the entries of the matrix ${\mathrm{M}}$ remains the same
for the matrix ${\mathrm{M}}^{-1}$. The commutation
relation $p_Ip_J=p_Jp_I+\varUpsilon$, $p_I\succ p_J$, $\varUpsilon\in\g$,
in $\U (\g )$ becomes
$\widetilde{\phantom{\,}p_Ip_J}=
\widetilde{\phantom{\,}p_Jp_I}+\widetilde{\phantom{\,}\varUpsilon\phantom{\,}}$,
$\widetilde{\phantom{\,}\varUpsilon\phantom{\,}}\in\p+\h$, in the
reduction algebra. Translate this into the ordering rule for the product $\mult$,
expressing the projections $\widetilde{\phantom{\,}pp\phantom{\,}}$'s in
terms of the products $\pp\mult\pp\ \, $'s with the help of the matrix  ${\mathrm{M}}^{-1}$
in both, left and right, hand sides: the right hand side, being rewritten in terms
of the multiplication $\mult$, consists of ordered terms only, the left hand side is
$\pp_I\mult \pp_J+\dots$, where dots stand for terms with $\pp_{I'}\mult \pp_{J'}$,
$p_I\succ p_{I'}$; such term is either ordered or, by induction in $I$, can be rewritten
in the ordered form as in \rf{not3}. The coefficient in front of
$\pp_{I'}\mult \pp_{J'}$ is from $\Uh$, so the reordering of the products
$\pp_{I'}\mult \pp_{J'}$ may force the coefficient of degree 1 or degree 0 term in \rf{not3}
to belong to $\Uh$.

\vskip .2cm
In the same manner we prove by induction on the filtration (described in the statement (b))
degree, that the algebra $\Zgk$  is generated over $\Uh$ by the elements $\{\pp_L\}$. To see this, consider the weight
basis, described in the statement (c), that is, the basis $\pp_{\overline{L}}$ ($\overline{L}$
is the multi-index) of the free $\Uh$-module $\Zgk^{(k)}/\Zgk^{(k-1)}$,
composed by images in $\Zgk$
of products $p_{L_1}^{n_1}p_{L_2}^{n_2}\cdots p_{L_m}^{n_m}$,
where $ p_{L_1}\prec p_{L_2}\prec\ldots\prec p_{L_m}$ and
$k=n_1+n_2+\dots +n_m$.
Equip the set of these basic elements with a total order $\prec$  compatible with the partial
order $<$ on $\h^*$; the compatibility has the same meaning as for
the elements $\{\pp_L\}$: [$\h$-weight of $\pp_{\overline{K}}$] $<$ [$\h$-weight of
 $\pp_{\overline{L}}$] $\Rightarrow$ $\pp_{\overline{K}}\prec
\pp_{\overline{L}}$. By the same, as above, arguments, referring to the structure of the
projector $P$, we have the following generalization of \rf{matrM}:
\be \pp_I\mult \pp_{\overline{J}}={\mathrm{M}}_{I\overline{J}K\overline{L}}
\widetilde{\phantom{\,}p_Kp_{\overline{L}}}\ ,\label{matrM2}\ee
where the  matrix ${\mathrm{M}}$, with entries in $\Uh$, has again
a triangular structure with 1's on the diagonal. Therefore, the matrix ${\mathrm{M}}$
is invertible and its inverse matrix ${\mathrm{M}}^{-1}$  has entries in $\Uh$.
The formula $\widetilde{\phantom{\,}p_Ip_{\overline{J}}}
={\mathrm{M}}_{I\overline{J}K\overline{L}}^{-1}\pp_K\mult \pp_{\overline{L}}$ implies the
induction step: the subspace $\Zgk^{(k+1)}$ is generated by products in $\Zgk$ of
elements from $\Zgk^{(1)}$. \hfill $\qed$

\vskip .2cm
Note that, before the localization, the algebra $\S(\g,\k)={\rm Norm}(\Ar\nplus)/\Ar\nplus$, as well as its image in $\Zgk$, is not
generated by the elements of degree 1. The subalgebra of $\S(\g,\k)$, generated by the elements of degree 1 (''step algebra''),
was the original subject of Mickelsson's investigation \cite{M}.

\begin{itemize}
\item[(e)] {The following monomials}  form a basis of the left $\Uh$-module $\Zgk$:
\be
 \pp_{I_1}\mult\pp_{I_2}\mult\cdots\mult\pp_{I_a},\qquad p_{I_1}\preceq p_{I_2}\preceq
 \ldots\preceq p_{I_a}\ .\label{inimpb}\ee
\end{itemize}

Before the proof of (e) we prove a more subtle statement.
\begin{prop}
Any expression in $\Zgk$ can be written in the ordered form by a repeated
application of \rf{not3} as instructions "replace the left hand side by the right hand side".
\end{prop}

{\it{Proof}} of Proposition. To save the space in the proof of this proposition we take a liberty to sometimes write
$I\prec J$ instead of $p_I\prec p_J$ (the same reservation concerns the use of $\preceq$, $\succ$ and $\succeq$).

\vskip .2cm
Consider the homogeneous quadratic part of the relations \rf{not3}:
\be\label{not3h} \pp_{I_1}\mult\pp_{I_2}=\displaystyle{\sum_{I_1',I_2':I_1'\preceq I_2'}}\dots\,
\pp_{I_1'}\mult\pp_{I_2'}\ ,\quad I_1\succ I_2\ ,\ee
where dots stand for coefficients from $\Uh$.
Denote by  ${\cal{I}}(\pp_{I_1}\mult\pp_{I_2})$
the right hand side of \rf{not3h}. We understand \rf{not3h} as
the set of instructions $\pp_{I_1}\mult\pp_{I_2}\leadsto {\cal{I}}(\pp_{I_1}\mult\pp_{I_2})$
($\leadsto$ stands for "replace") in the free algebra with
the weight generators $\pp_I$.

Let us prove the statement for a cubic monomial $\pp_{I_1}\mult\pp_{I_2}\mult\pp_{I_3}$.
For such a monomial one can apply the instructions \rf{not3h} to
$\pp_{I_1}\mult\pp_{I_2}$ if $I_1\succ  I_2$ and to $\pp_{I_2}\mult\pp_{I_3}$ if $I_2\succ  I_3$.
Denote the results by
${\cal{I}}_{12}(\pp_{I_1}\mult\pp_{I_2}\mult\pp_{I_3})$ and
${\cal{I}}_{23}(\pp_{I_1}\mult\pp_{I_2}\mult\pp_{I_3})$ respectively.

For an element $\psi\in\h^*$, $\psi =\sum l_i\alpha_i$, where $\alpha_i$
are the simple roots, let $d(\psi ):=\sum l_i$. The function $d$ is
compatible with the partial order $<$ on $\h^*$ in the sense that $d(\alpha )<d(\beta )$ if
 $\alpha <\beta$. Denote by the same letter $d$ the
function on the set of indices, labeling the weight base of $\p$; it is defined by
$d(I):=d(\mu_I)$, where $\mu_I$ is the weight of $\pp_I$.

We have
$d(I_1')+d(I_2')=d(I_1)+d(I_2)$ for any monomial $\pp_{I_1'}\mult\pp_{I_2'}$ appearing in
the right hand side of \rf{not3h} (and the difference
$d(I_1)-d(I_1')$ is an integer). Since $I_1\succ I_2$ and $I_1'\preceq I_2'$, it follows
that $d(I_1)\geq d(I_2)$ and $d(I_1')\leq d(I_2')$; therefore, $d(I_1')\leq d(I_1)$
and $d(I_2')\geq d(I_2)$.

Associate to a monomial
 $\pp_{I_1}\mult\pp_{I_2}\mult\pp_{I_3}$, that is,
to an ordered triple $(I_1,I_2,I_3)$ of indices, the number
${\mathfrak{d}}(I_1,I_2,I_3):=2d(I_1)+d(I_2)$.
When we apply the ordering instructions ${\cal{I}}_{12}$ or ${\cal{I}}_{23}$ to
$\pp_{I_1}\mult\pp_{I_2}\mult\pp_{I_3}$, the
function ${\mathfrak{d}}$ does not increase; that is, the value of ${\mathfrak{d}}$
on any of the appearing monomials is not greater than
${\mathfrak{d}}(I_1,I_2,I_3)$. Indeed, if we replace $\pp_{I_1}\mult\pp_{I_2}$ by
 $\pp_{I_1'}\mult\pp_{I_2'}$ then $2d(I_1')+d(I_2')=d(I_1')+
\bigl( d(I_1')+d(I_2')\bigr)=d(I_1')+\bigl( d(I_1)+d(I_2)\bigr)\leq d(I_1)+
\bigl( d(I_1)+d(I_2)\bigr)=2d(I_1)+d(I_2)$; and if we replace $\pp_{I_2}\mult\pp_{I_3}$
by $\pp_{I_2'}\mult\pp_{I_3'}$ then simply $d(I_2')\leq d(I_2)$ and $d(I_1')=d(I_1)$.

 For a linear combination $X=\sum c_{I_1I_2I_3}\pp_{I_1}\mult\pp_{I_2}\mult\pp_{I_3}$
 of cubic monomials, with coefficients $c_{I_1I_2I_3}\in\Uh$, denote
the maximal value of ${\mathfrak{d}}$ on the monomials, appearing in $X$, by the same symbol
 ${\mathfrak{d}}$; that is, ${\mathfrak{d}}(X):={\displaystyle{\max_{(I_1,I_2,I_3):
 c_{I_1I_2I_3}\neq 0}}}{\mathfrak{d}}(I_1,I_2,I_3)$.

\vskip .2cm Assume that the assertion is false and there exists a
cubic monomial which cannot be ordered by the instructions
\rf{not3h}. Since $\k +\p$ is finite-dimensional, the set of values
of the function ${\mathfrak{d}}$ on cubic monomials is bounded from
below. So the minimal value ${\mathfrak{d}}_{\min}$ of the function
${\mathfrak{d}}$ on the set of cubic monomials which cannot be
ordered is finite, ${\mathfrak{d}}_{\min} >-\infty$. Let
$\pp_{I_1}\mult\pp_{I_2}\mult\pp_{I_3}$ be a monomial, which cannot
be ordered, with
${\mathfrak{d}}(I_1,I_2,I_3)={\mathfrak{d}}_{\min}$. The application
of the ordering instructions \rf{not3h} cannot strictly decrease the
value of ${\mathfrak{d}}$, this would contradict to the minimality
of ${\mathfrak{d}}(I_1,I_2,I_3)$. Therefore, among the appearing
monomials, there is at least one monomial
$\pp_{I_1'}\mult\pp_{I_2'}\mult\pp_{I_3'}$ with the same value of
${\mathfrak{d}}$. If $\pp_{I_1'}\mult\pp_{I_2'}\mult\pp_{I_3'}$
appears in ${\cal{I}}_{12}(\pp_{I_1}\mult\pp_{I_2}\mult\pp_{I_3})$
then $2d(I_1')+d(I_2')=2d(I_1)+d(I_2)$,  ${ {d(I_1)\geq d(I_2),\ d(I_1')\leq d(I_2')}}$
 and $I_3'=I_3$; since the
total weight is conserved, $d(I_1')=d(I_1)$, { and} ${{d(I_2')=d(I_2)}}$, { so}
$d(I_1)=d(I_2)=d(I_1')=d(I_2')$. If
$\pp_{I_1'}\mult\pp_{I_2'}\mult\pp_{I_3'}$ appears in
${\cal{I}}_{23}(\pp_{I_1}\mult\pp_{I_2}\mult\pp_{I_3})$ then
${ {d(I_2)\geq d(I_3),\ d(I_2')\leq d(I_3')}}$, $d(I_2')=d(I_2)$ and $I_1'=I_1$;
 { by the same arguments we have again},
$d(I_2)=d(I_3)=d(I_2')=d(I_3')$. Due to the structure of the matrix
${\mathrm{M}}$, defined in \rf{matrM}, and the arguments used in the
proof of the statement (d), ${\cal{I}}(\pp_{I}\mult\pp_{J})$ with
$d(I)=d(J)$ contains exactly one monomial $\pp_{I'}\mult\pp_{J'}$
with $d(I')=d(I)$ and this monomial is $\pp_{J}\mult\pp_{I}$.
Therefore, up to monomials with the value of ${\mathfrak{d}}$
smaller than ${\mathfrak{d}}_{\min}$ (they can be ordered by
assumption) and up to a coefficient from $\Uh$, the operation
${\cal{I}}_{12}$, $I_1\succ I_2$, is simply
$\pp_{I_1}\mult\pp_{I_2}\leadsto\pp_{I_2}\mult\pp_{I_1}$; the
operation ${\cal{I}}_{23}$, $I_2\succ I_3$, is
$\pp_{I_2}\mult\pp_{I_3}\leadsto\pp_{I_3}\mult\pp_{I_2}$.
 The transpositions (12) and (23) of neighbors generate all permutations of three letters.
The orbit of
$\pp_{I_1}\mult\pp_{I_2}\mult\pp_{I_3}$ under the group of permutations of three letters
$I_1,I_2$ and $I_3$ contains the ordered monomial,
the contradiction.

The degree 0 or 1 terms, contained in the full instructions \rf{not3}, may only cause
an appearance of linear or quadratic terms in the process of ordering of a cubic
polynomial. So, any cubic polynomial can be ordered by \rf{not3} as well.

More generally, to a monomial $X=\pp_{I_1}\mult\pp_{I_2}\mult\cdots\mult\pp_{I_k}$
of an arbitrary degree $k$ we  associate  the number
${\mathfrak{d}}(I_1,\ldots,I_k):=(k-1)d(I_1)+(k-2)d(I_2)+\ldots+d(I_{k-1})$,
and, in the minimal situation, conclude that
up to terms smaller than $X$ in an appropriate sense, the instructions
essentially reduce to transpositions $(i,i+1)$ of neighbors, which generate the whole
symmetric group on $k$ letters, and thus an ordered expression is in the orbit.
\hfill $\qed$

\vskip .2cm {\it{Proof}} of statement (e).   By the statement (d)
above, the algebra $\Zgk$ is
 generated by $\pp_I$
and, due to the form \rf{not3} of relations, has a filtration by the $\mult$-degree.
Let $\Zgk^{(\mult k)}$
be the subspace of elements of degree not greater
than $k$ with respect to the product $\mult$. Since
$\pp_{I_1}\mult\pp_{I_2}\mult\dots\mult\pp_{I_k}=\pp_{I_1}P\pp_{I_2}P\dots P\pp_{I_k}$,
it follows that
$\Zgk^{(\mult k)}\subset \Zgk^{(k)}$. The opposite inclusion holds as well because the algebra
$\Zgk$ is generated by $\pp_I$. We conclude that the two
filtrations coincide.

\vskip .2cm
Therefore, every element $\widetilde{\phantom{\,}p_{I_1}\cdots p_{I_k}}$,
$I_1\preceq \ldots\preceq I_k$, is in $\Zgk^{(\mult k)}$
and, by proposition above, can be ordered.
The cardinalities of the sets  $\{\widetilde{p_{I_1}\cdots p_{I_k}}\mid I_1\preceq \
\ldots \preceq I_k\}$ and
$\{\pp_{I_1}\mult\ldots\mult\pp_{I_k}\mid I_1\preceq \ldots\preceq I_k\}$ are equal,
so due to \rf{inimp} the set $\{\pp_{I_1}\mult\ldots\mult\pp_{I_k}\mid
I_1\preceq \ldots\preceq I_k\}$ is a basis of
$\Zgk^{(\mult k)}/\Zgk^{(\mult (k-1))}$.\hfill $\qed$

\vskip .2cm
Note that for an order which is not compatible with the partial order $<$ on $\h^*$,
the ordering relations of the form \rf{not3} may exist but the statement (e) does not
necessarily hold. For instance, the ordering relations \rf{not3} can be written for a
lexicographical order for the generators $\z_{ij}$ and $t_i$ (with $\z_{ii}=t_i$) of the algebra
$\Z_n$, defined in the next Section, but the ordering procedure loops for cubic monomials,
already for $n=2$ (we don't give details; it is an explicit calculation).

\setcounter{equation}{0}
\section{Diagonal reduction algebra of $\gl_n$\vspace{.25cm}}\label{section-notation}
\vskip -.2cm
Let $\gl_n$ be the Lie algebra of the general linear group of $n$-dimensional complex linear space.
  Consider the reductive pair $(\g,\k)$ with $\g=\gl_n\oplus\gl_n$ and $\k=\gl_n$  diagonally
embedded into $\gl_n\oplus\gl_n$. The corresponding reduction algebra we call  'diagonal
reduction algebra'
and denote it by $\Z_n$.

We fix the following notations for generators of these algebras $\g$
and $\k$. Let $E_{ij}^{(1)}$ and $E_{ij}^{(2)}$, $i,j=1\lcd n$, be
the standard generators of the two copies of the Lie algebra $\gl_n$
in $\gl_n\oplus\gl_n$,
$$[E_{ij}^{(a)},E_{kl}^{(b)}]=\delta_{ab}\left(\delta_{jk}E_{il}^{(a)}-
\delta_{il}E_{kj}^{(a)}\right)\ ,$$
where $\delta_{ab}$ and $\delta_{ij}$ are the Kronecker symbols. Set
\be e_{ij}:=\frac{1}{2}(E_{ij}^{(1)}+E_{ij}^{(2)})\ ,\qquad E_{ij}:=\frac{1}{2}(E_{ij}^{(1)}-
E_{ij}^{(2)})\ .\ee
The elements $e_{ij}$ span the diagonally embedded Lie algebra $\k\simeq\gl_n$,
 while $E_{ij}$ form an adjoint $\k$-module $\mathfrak{p}$. The
Lie algebra $\k$ and the space $\mathfrak{p}$ constitute a symmetric pair,
that is, $[\k,\k]\subset \k$, $[\k,\p]\subset \p$, and $[\p,\p]\subset \k$:
\be\begin{array}{ll}[e_{ij},e_{kl}]&=\delta_{jk}e_{il}-\delta_{il}e_{kj}\ ,
\qquad [e_{ij},E_{kl}]=\delta_{jk}E_{il}-\delta_{il}E_{kj}\ ,\\[.5em]
 [E_{ij},E_{kl}]&=\delta_{jk}e_{il}-\delta_{il}e_{kj}\ .\end{array}\ee
In the sequel, $h_a$ means the element $e_{aa}$ of the Cartan subalgebra $\h$ of the
subalgebra $\k\in\gl_n\oplus\gl_n$ and $h_{ab}$ the
element $e_{aa}-e_{bb}$.

\vskip .2cm
Let $\{\ve_a\}$ be the basis of $\h^*$ dual to the basis $\{ h_a\}$ of $\h$,
$\ve_a(h_b)=\delta_{ab}$. We shall use as well the root notation $h_\alpha$,
$e_{\alpha}$, $e_{-\alpha}$ for elements of $\k$, and $H_\alpha$, $E_{\alpha}$, $E_{-\alpha}$
for elements of $\p$. The Lie sub-algebra $\nplus$ in the
triangular decomposition \rf{intro1} is spanned by the root vectors $e_{ij}$ with $i<j$
and the Lie sub-algebra $\nminus$ by the root vectors $e_{ij}$ with $i>j$. Let $\bb_+$ and $\bb_-$
be the corresponding Borel sub-algebras, $\bb_+=\h\oplus\nplus$, $\bb_-=\h\oplus\nminus$.
The system  $\Delta_+$ of positive roots of $\k$ consists of roots $\ve_i-\ve_j$ with $i<j$
and the system
$\Delta_-$ consists of roots
$\ve_i-\ve_j$ with $i>j$.

We fix the following action of the cover of the symmetric group $\S_n$ (the Weyl group of the
diagonal $\k$) on the Lie algebra $\gl_n\oplus\gl_n$ by automorphisms
\be\sy_i(x):=\mathrm{Ad}_{\exp(e_{i,i+1})}\mathrm{Ad}_{\exp(-e_{i+1,i})}\mathrm{Ad}_{\exp(e_{i,i+1})}(x)\ ,\ee
so that
 $\sy_i(e_{kl})=(-1)^{\delta_{ik}+\delta_{il}}e_{\sigma_i(k)\sigma_i(l)}$ and $
 \sy_i(E_{kl})=(-1)^{\delta_{ik}+\delta_{il}}E_{\sigma_i(k)\sigma_i(l)}$.
Here $\sigma_i=(i,i+1)$ is an elementary transposition in the symmetric group.
 We extend naturally the above action of the cover of $\S_n$ to the action by
automorphisms on the associative algebra $\Ar\equiv\Ar_n:=\U(\gl_n)\otimes \U(\gl_n)$.
The restriction of this action to $\h$ coincides
with the natural action $\si(h_{k})=h_{\si(k)}$, $\si\in\S_n$,  of the Weyl group on the
Cartan sub-algebra.
The shifted action \rf{shifted} of the Weyl group on $\h$ looks as:
\be\label{not1}\si\circ h_k:= h_{\si(k)}+k-\si(k)\ ,\qquad k=1,...,n\ ;\quad \si\in\S_n\ .\ee
It becomes the usual action for the variables
\be \th_k:=h_k-k\ ,\qquad \th_{ij}:=\th_i-\th_j\ ;\label{hrond}\ee
so that for any $\si\in\S_n$ we have
$\si\circ\th_k=\th_{\si(k)}$ and\ $\si\circ\th_{ij}=\th_{\si(i)\si(j)}$ .
The set of denominators, defining the localizations $\Uh$ and $\Ab$ consists of elements
\be h_{ij}+l\ ,\qquad l\in\ZZ\ , \quad 1\leq i<j\leq n\ .\label{musel}\ee
 We choose the set of vectors $E_{ij}$, $i,j=1,...,n$, as a basis of the space $\p$.
 The weight of $E_{ij}$ is $\ve_i-\ve_j$.
The compatibility of a total order $\prec$ with the partial order $<$ on $\h^*$
means the condition
\begin{equation}\label{not4a}E_{ij}\prec E_{kl}\qquad \text{if}\qquad i-j>k-l\ .\end{equation}
The order in each subset $\{E_{ij}|i-j=a\}$ with a fixed $a$ can be chosen arbitrarily.
For instance, we can set
\begin{equation}\label{not4}E_{ij}\prec E_{kl}\quad \text{if}\quad i-j>k-l\quad
\text{or}\quad i-j=k-l\quad\text{and}\quad i>k\ .\end{equation}
\vskip .2cm
 Denote the images of the elements $E_{ij}$ in $\Z_n$ by $\s_{ij}$.
  We use also the notation $t_i$ for the elements $\s_{ii}$ and $t_{ij}:=t_i-t_j$
  for the elements
   $\s_{ii}-\s_{jj}$.
The order \rf{not4} induces as well the order on the generators $\s_{ij}$ of
the algebra $\Z_n$:
\be\label{not4b}\s_{ij}\prec\s_{kl}\ \Leftrightarrow\ E_{ij}\prec E_{kl}\ .\ee
The  statement \rf{not3} implies an existence of structure constants
$\mathrm{B}_{(ab),(cd),(ij),(kl)} \in \Uh$ and $\mathrm{D}_{(ab),(cd)}\in\Uh$
such that for any $a,b,c,d=1,\ldots,n$ we have
\be\label{not6}\s_{ab}\mult\s_{cd}=\sum_{i,j,k,l:\s_{ij}
\preceq\s_{kl}}\mathrm{B}_{(ab),(cd),(ij),(kl)}\s_{ij}\mult\s_{kl}+\mathrm{D}_{(ab),(cd)}\ .\ee
Linear terms in the right hand side of \rf{not6} are absent since here
$(\g,\k)$ form a symmetric pair.
The relations \rf{not6} together with the weight conditions
\be [h,\s_{ab}]=(\ve_a-\ve_b)(h) \s_{ab}\label{wede}\ee
are the defining relations for the algebra $\Z_n$.

\vskip .2cm The structure of denominators of entries of the matrices
${\mathrm{M}}$ and ${\mathrm{M}}^{-1}$, mentioned in the proof of
\rf{not3} above, shows that for the algebra $\Z_n$ the denominators
of the structure constants $\mathrm{B}_{(ab),(cd),(ij),(kl)}$ and
$\mathrm{D}_{(ab),(cd)}$ are products of linear factors of the form
$\th_{ij}+\ell$, $i<j$, where $\ell\geq -1$ is an integer. This is
because in our situation the $\sl_2$--sub-algebra (of the diagonal
$\gl_n$), corresponding to an arbitrary positive root $\ve_i-\ve_j$,
$i<j$, has only  1, 2- and 3-dimensional representations in
$\p$, so the numbers $\ell$'s in the denominators of the summands of
the projector can drop at most by 2 due to the presence of the term
$(\pi ,\gamma )$.

 The Chevalley anti-involution $\epsilon$ in $\U(\gl_n\oplus\gl_n)$,
$\epsilon(e_{ij}):=e_{ji}$, $\epsilon(E_{ij}):=E_{ji}$, induces the anti-involution $\epsilon$ in
the algebra $\Z_n$:
\be\epsilon(\s_{ij})=\s_{ji}\ ,\qquad \epsilon(h_k)=h_k\ .\label{anep}\ee
Besides, the outer automorphism of the Dynkin diagram of $\gl_n$ induces the involutive automorphism
$\omega$ of $\Z_n$,
\be\label{not2a}\omega(\s_{ij})=(-1)^{i+j+1}\s_{j'i'}\ ,\qquad \omega (h_k)=-h_{k'}\ ,\ee
where $i'=n+1-i$. The operations $\epsilon$ and $\omega$ commute, $\epsilon\omega =\omega\epsilon$.

\vskip .2cm
Central elements of the sub-algebra $\U(\gl_n)\otimes 1\subset \Ar$, generated by $n$ Casimir
operators of degrees $1\lcd n$,  as well as central elements of the sub-algebra
$1\otimes \U(gl_n)\subset \Ar$ project to central elements of the algebra $\Z_n$. In particular,
central elements of degree $1$ project to central elements
\be h_1+\ldots +h_n\qquad
\text{and}\qquad
\label{clit}t_1+\ldots +t_n\ee
of the algebra $\Z_n$. The difference of central elements of degree two projects to the central
element
\be\label{drclit}\sum_{i=1}^n(h_i-2i)t_i\ee
of the algebra $\Z_n$. The images of other Casimir operators are more complicated.

\subsection{Change of variables}

We shall use the following elements of $\Uh$:
$$\!\! \tphi_{ij}:=\frac{\th_{ij}}{\th_{ij}-1}\ ,\ \tpsi_{ij}:=\frac{\th_{ij}-1}{\th_{ij}}\ ,
\ \ttphi_{ij}:=\frac{\th_{ij}-1}{\th_{ij}-2}\ ,\ \ttpsi_{ij}:=
\frac{\th_{ij}-2}{\th_{ij}-1}\ ,\ \Cprime_{ij}:=\frac{\th_{ij}-3}{\th_{ij}-2}\ ,$$
the variables $\th_{ij}$ are defined in \rf{hrond}. Note that $\ \tphi_{ij}\tpsi_{ij}=\ttphi_{ij}\ttpsi_{ij}=1$.

\vskip .2cm
Define  elements $\tt_1,\ldots, \tt_n\in\Z_n$ by
\be\label{3.1}\tt_1:=t_1\ ,\quad\tt_2:=\q_{1}(t_1)\ ,\quad\tt_3:=\q_{2}\q_{1}(t_1)\ ,
\quad\ldots\quad,\tt_n:=\q_{n-1}\cdots \q_{2}\q_{1}(t_1)\ .\ee
Using \rf{not7} we find the relations
\be\begin{array}{ll}&\displaystyle{\q_i(t_i)\ =-\frac{1}{\th_{i,i+1}-1}t_i+
\frac{\th_{i,i+1}}{\th_{i,i+1}-1}t_{i+1}\ ,}\\[1.3em]
&\displaystyle{\q_i(t_{i+1})=\ \frac{\th_{i,i+1}}{\th_{i,i+1}-1}t_i-\frac{1}{\th_{i,i+1}-1}
t_{i+1}\ ,}\\[1.8em]
&\displaystyle{\q_i(t_k)\ =\quad t_k\ ,\qquad\qquad k\not=i,i+1}\ ,\end{array}\label{acwot}\ee
which can be used to convert the definition \rf{3.1} into  a linear over the ring $\Uh$ change
of variables:
\be\begin{array}{rl}\label{3.1a}\ds{\tt_l}&=\ds{t_l\prod_{j=1}^{l-1}\tphi_{jl}
-\sum_{k=1}^{l-1}t_k\frac{1}{\th_{kl}-1}\prod_{j=1}^{k-1}\tphi_{jl}\ ,}\\[2.2em]
\ds{t_l}&=\ds{\tt_l\prod_{j=1}^{l-1}\tpsi_{jl} +\sum_{k=1}^{l-1}\tt_k\frac{1}{\th_{kl}}
\prod_{\begin{array}{c}\scriptstyle{j=1}\\[-.5em] \scriptstyle{j\neq k}\end{array}}^{l-1}
\tpsi_{jk}\ .}
\end{array}\ee

In terms of the new variables $\tt$'s, the linear in $t$ central element (\ref{clit}) reads
$$\sum t_i=\sum\tt_i\prod_{a:a\neq i}\frac{\th_{ia}+1}{\th_{ia}}\ .$$

In the following, we use the notion of {\it coefficient-bounded} formulas and relations.
It means the following.
 Given a family of formulas for each $n$ (expressing some
action, relations {\it etc.}) with coefficients in $\Uh$, we say that it is
 coefficient-bounded if the degrees of the numerators and denominators (in
the reduced form, with no common factors) of the coefficients do not grow with
$n$.

 For example, the set of relations for  $\Z_n$, which we shall
exhibit, will have coefficient-bounded terms with respect to a certain set of
generators.
In this sense the action \rf{acwot} is coefficient-bounded while the change of variables
\rf{3.1a} is however not
coefficient-bounded.
\subsection{Braid group action}\label{brgac}

 Since
$\q_{i}^2(x)=x$ for any element $x$ of zero weight, the braid group acts as its symmetric group
quotient on the space of weight 0 elements.
Although the change of variables \rf{3.1a} is not coefficient-bounded in the sense of
Section \ref{section2new}, the action of the
transformations $\q_{i}$ on the new variables $\tt$'s is coefficient-bounded: it follows from
\rf{3.1} and $\q_{i}(t_1)=t_1$ for all $i>1$ that
\be\label{3.2}\q_{\sigma}(\tt_i)=\tt_{\sigma(i)}\qquad
\text{for any}\qquad \sigma\in \S_n. \ee

The action of the Zhelobenko automorphisms  on the generators $\s_{kl}$
 looks as follows:
\begin{align}&\q_i(\s_{ik})=-\s_{i+1,k}\tphi_{i,i+1}\ ,&& \notag
\q_i(\s_{ki})=-\s_{k,i+1}\ ,&&k\not=i,i+1\ ,\\ \label{3.5}
&\q_i(\s_{i+1,k})=\s_{i,k}\ ,&&\q_i(\s_{k,i+1})=\s_{k,i}\tphi_{i,i+1}\ ,&&k\not=i,i+1\ ,\\
&\q_i(\s_{i,i+1})=-\s_{i+1,i}\tphi_{i,i+1}\ttphi_{i,i+1}\ ,&& \notag
\q_i(\s_{i+1,i})=-\s_{i,i+1}\ ,\\ \notag
&\q_{i}(\s_{j,k})=\s_{j,k}\ ,\ \ j,k\neq i,i+1\ .&&\end{align}

\vskip .1cm
Denote $i'=n+1-i$, as before. The braid group action \rf{3.5} is compatible with the anti-involution
$\epsilon$
and the involution $\omega$ (note that $\omega (\th_{ij})=\th_{j'i'}$),  see \rf{anep} and \rf{not2a}, in the following sense:
\begin{align} \epsilon\, \q_i&=\q_i^{-1}\epsilon\ \ ,\label{qeps} &\omega  \q_i&=\q_{i'-1}\omega\ .
\end{align}

\vskip .1cm
Let $w_0$ be the longest element of the Weyl group of $\gl_n$, the symmetric group $\S_n$.
Similarly to the squares of
the transformations corresponding to the simple roots, see \rf{invr}, the action of $\q_{w_0}^2$
is the conjugation by a certain element of $\Uh$. Moreover, one can observe by a direct calculation,
that

\be \label{3.6}\q_{w_0}(\s_{ij})=(-1)^{i+j} \s_{i'j'}\prod_{a:a<i'}\tphi_{ai'}\
\prod_{b:b>j'}\tphi_{j'b}\ ,\qquad
\q_{w_0}(\tt_i)=\tt_{i'}\ .
\ee

The formula \rf{3.6} implies the existence of the ordering relations for the generators $\s_{ij}$
in the inverse to \rf{not4}-\rf{not4b} order.
\begin{cor}\hspace{-.2cm}. There exist $\mathrm{B}'_{(ab),(cd),(ij),(kl)}$ and
$\mathrm{D}'_{(ab),(cd)}\in\Uh$
such that for any $\s_{ab}$ and $\s_{cd}$ we have
\begin{equation}\label{not6a}\s_{ab}\mult\s_{cd}=\sum_{i,j,k,l:\s_{kl}
\preceq\s_{ij}}\mathrm{B}'_{(ab),(cd),(ij),(kl)}
\s_{ij}\mult\s_{kl}+\mathrm{D}'_{(ab),(cd)}\ .\end{equation}
\label{opor}\end{cor}

Indeed, we apply the transformation $\q_{w_0}$ to the equalities \rf{not6} and substitute \rf{3.6}.
This gives the relations \rf{not6a} since the assignment $(i,j)\mapsto(i',j')$ reverses the
order $\prec$.

\subsection{Defining relations}\label{section3.3}

To save  space we omit in this section the symbol $\mult$ for the multiplication in the algebra
$\Z_n$. It should not lead to any confusion since no other multiplication is used in this section.

\vskip .2cm Each relation which we will derive will be of a certain
weight, equal to a sum of two roots. From general considerations the
upper estimate for the number of terms in a quadratic relation of
weight $\lambda=\alpha+\beta$ is the number $|\lambda|$ of quadratic
combinations $\s_{\alpha'}\s_{\beta'}$ with
$\alpha'+\beta'=\lambda$. There are several types of relation
weights, excluding the trivial one, $\lambda=2(\ve_i-\ve_j)$,
$|\lambda|=1$:
\begin{enumerate}
\item $\lambda=\pm(2\ve_i-\ve_j-\ve_k)$, where $i,j$ and $k$ are pairwise distinct.
Then $|\lambda|=2$.
\item $\lambda=\ve_i-\ve_j+\ve_k-\ve_l$ with pairwise distinct
$i,j,k$ and $l$. Then $|\lambda|=4$.
\item $\lambda=\ve_i-\ve_j$, $i\neq j$. For $\s_{\alpha'}\s_{\beta'}$, there are $2(n-2)$
possibilities (subtype 3a) with $\alpha'=\ve_i-\ve_k$,
$\beta'=\ve_k-\ve_j$ or $\alpha'=\ve_k-\ve_j$, $\beta'=\ve_i-\ve_k$ with $k\neq i,j$ and $2n$
possibilities (subtype 3b) with $\alpha'=0$,
$\beta'=\ve_i-\ve_j$ or $\alpha'=\ve_i-\ve_j$, $\beta'=0$. Thus $|\lambda|=4(n-1)$.
\item $\lambda=0$. There are $n^2$ possibilities (subtype 4a) with $\alpha'=0$, $\beta'=0$
and $n(n-1)$ possibilities (subtype 4b) with
$\alpha'=\ve_i-\ve_j$,  $\beta'=\ve_j-\ve_i$, $i\neq j$. Here $|\lambda|=n(2n-1)$.
\end{enumerate}

Below we write down relations for each type (and subtype) separately. The relations of  types 1
and 2 have a simple form in terms of the original
generators $\s_{ij}$. To write the relations  of  types 3 and 4, it is convenient to renormalize
the generators $\s_{ij}$
with $i\not=j$. Namely, we set
\begin{equation}\label{not8}\hs_{ij}=\s_{ij}\prod_{k=1}^{i-1}\tphi_{ki}\ .\end{equation}

In terms of the generators $\hs_{ij}$, the formulas \rf{3.5} for the action of the automorphisms
$\q_i$  translate as follows:
\begin{equation*}\begin{array}{lll}\q_{i}(\hs_{ik})=-\hs_{i+1,k}\ ,&
\q_{i}(\hs_{i+1,k})=\hs_{i,k}\tphi_{i+1,i}\ ,&k\not=i,i+1\ ,\\[.3 em]
\q_{i}(\hs_{ki})=-\hs_{k,i+1}\ ,&
\q_{i}(\hs_{k,i+1})=\hs_{k,i}\tphi_{i,i+1}=\tpsi_{i+1,i}\hs_{k,i}\ ,&k\not=i,i+1\ ,\\[.3 em]
\q_{i}(\hs_{i,i+1})=-\tpsi_{i+1,i}\hs_{i+1,i}\ ,&\q_{i}(\hs_{i+1,i})=-\hs_{i,i+1}\tphi_{i+1,i}\
,&\\[.3 em]
\q_{i}(\hs_{j,k})=\hs_{j,k}\ , \ \ j,k\neq i,i+1\ .&&\end{array}\end{equation*}
Although the renormalization \rf{not8} is not coefficient-bounded, the action of the braid group
stays coefficient-bounded.

\paragraph{1.} The relations of  the type 1 are:
\begin{equation}\label{relation1}\s_{ij}\s_{ik}=\s_{ik}\s_{ij}\tphi_{kj}\ ,\qquad
\s_{ji}\s_{ki}=\s_{ki}\s_{ji}\tpsi_{kj}\ ,\qquad \text{for}\quad j<k\ ,\ i\not=j,k\ .\end{equation}

\paragraph{2.} Denote
$D_{ijkl}:= \th_{ik}^{-1}-\th_{jl}^{-1}$.
Then, for any four pairwise different indices $i,j,k$ and $l$, we have the following relations of
the type 2:
\begin{equation}\label{relation2}\begin{split} [ \s_{ij},\s_{kl} ]&=\s_{kj}\s_{il}D_{ijkl}\ ,
\qquad i<k\ ,\ j<l \ ,\\[.4em] \s_{ij}\s_{kl}-\s_{kl}\s_{ij}\tphi_{jl}'\tphi_{lj}'&=\s_{kj}\s_{il}
D_{ijkl}\ ,\qquad i<k\ ,\ j>l\ .\end{split}\end{equation}

\paragraph{3a.} Let $i\neq k\neq l\neq i$. Denote
$$\mathring{E}_{ikl}:=-\left((\tt_i-\tt_k)\frac{\th_{il}+1}{\th_{ik}\th_{il}}+(\tt_k-\tt_l)
\frac{\th_{il}-1}{\th_{kl}\th_{il}}\right)\hs_{il}+\sum_{a:a\neq i,k,l}
\hs_{al}\hs_{ia}\frac{\ttphi_{ai}}{\th_{ka}+1}\ .$$
With this notation the first group of the relations of the type 3 is:
\begin{align}\notag\hs_{ik}\hs_{kl}\tpsi_{ik} -\hs_{kl}\hs_{ik}\ttphi_{ki}
&=\mathring{E}_{ikl}\ , && i<k<l\ ,\\[.4em]  \notag
\hs_{ik}\hs_{kl}\tpsi_{ik}\tpsi_{lk}\ttphi_{lk} -\hs_{kl}\hs_{ik}\ttphi_{ki}
&=\mathring{E}_{ikl}\ , && i<l<k\ ,\\[.4em]  \label{relation3}
\hs_{ik}\hs_{kl}\tphi_{ki} -\hs_{kl}\hs_{ik}\ttphi_{ki}&=\mathring{E}_{ikl}\ , && k<i<l\ ,\\[.4em]
\notag
\hs_{ik}\hs_{kl}\tphi_{ki}\tphi_{li}\ttpsi_{li} -\hs_{kl}\hs_{ik}\ttphi_{ki}
&=\mathring{E}_{ikl}\ ,&& k<l<i\ ,\\[.4em]  \notag
\hs_{ik}\hs_{kl}\tpsi_{ik}\tpsi_{lk}\ttphi_{lk}\tphi_{li}\ttpsi_{li}-\hs_{kl}\hs_{ik}\ttphi_{ki}
&=\mathring{E}_{ikl}\ , && l<i<k\ ,\\[.4em]  \notag
\hs_{ik}\hs_{kl}\tphi_{ki}\tpsi_{lk}\ttphi_{lk}\tphi_{li}\ttpsi_{li}
-\hs_{kl}\hs_{ik}\ttphi_{ki}&=\mathring{E}_{ikl}\ , && l<k<i\ .\end{align}

The relations \rf{relation3} can be written in a more compact way with the help of both systems,
$\s_{ij}$ and $\hs_{ij}$, of generators. Let now
$$E_{ikl}:=-\left((\tt_i-\tt_k)\frac{\th_{il}+1}{\th_{ik}\th_{il}}+(\tt_k-\tt_l)
\frac{\th_{il}-1}{\th_{kl}\th_{il}}\right)\s_{il}+\sum_{a:a\neq i,k,l}
\hs_{al}\s_{ia}\frac{\ttphi_{ai}}{\th_{ka}+1}\ .$$
Then
\be\begin{array}{rcll} \s_{ik}\hs_{kl}\tpsi_{ik} -\hs_{kl}\s_{ik}\ttphi_{ki}&=&E_{ikl}\ ,\ \ & k<l\ ,
\\[.4em]
\s_{ik}\hs_{kl}\tpsi_{ik}\tpsi_{lk}\ttphi_{lk} -\hs_{kl}\s_{ik}\ttphi_{ki}&=&E_{ikl}\ ,\ \ & l<k\ .
\end{array}\label{shfo3a}\ee
Moreover, after an extra redefinition:
$\!\mathring{\hspace{.1cm}z}_{kl}\!\!\!\!\!\!\mathring{\phantom{z}}\ \ =\hs_{kl}
\ttphi_{lk}$ for $k>l$, the left hand side of the
second line in \rf{shfo3a} becomes, up to a common factor, the same as the left hand side
of the first line, namely, it reads $(\s_{ik}
\mathring{\hspace{.1cm}z}_{kl}\!\!\!\!\!\!\mathring{\phantom{z}}\ \ \tpsi_{ik} -
\mathring{\hspace{.1cm}z}_{kl}\!\!\!\!\!\!\mathring{\phantom{z}}\ \ \s_{ik}\ttphi_{ki})
\tpsi_{lk}$.

\paragraph{3b.} Let $l\neq j$. The second group of relations of the type 3 reads:
\begin{align}\notag\hs_{ij}\tt_{i}=&\ \tt_i\hs_{ij}\Cprime_{ji}-\tt_j\hs_{ij}
\frac{1}{\th_{ij}+2}-
\sum_{a:a\not=i,j}\hs_{aj}\hs_{ia}\frac{1}{\th_{ia}+2}\ ,\\[.4em]
\label{3.12}\hs_{ij}\tt_{j}=&-\tt_i\hs_{ij}\frac{\Cprime_{ji}}{\th_{ij}-1}+
\tt_j\hs_{ij}\tphi_{ij}\tpsi_{ji}\ttphi_{ji}+
\sum_{a:a\not=i,j}\hs_{aj}\hs_{ia}\tphi_{ij}\tpsi_{ji}\frac{\ttphi_{ai}}{\th_{ja}+1}\ ,\\[.4em]
\notag
\hs_{ij}\tt_{k}=&\,\tt_i\hs_{ij}\frac{(\th_{ij}+3)\ttphi_{ji}}{(\th_{ik}^2-1)(\th_{jk}-1)}+
\tt_j\hs_{ij}\frac{(\th_{ij}+1)\ttphi_{ji}}{(\th_{ik}-1)(\th_{jk}-1)^2}+
\tt_k\hs_{ij}\tphi_{ik}\tphi_{ki}\tphi_{jk}\ttpsi_{jk}\\[.4em]  \notag
-&\hs_{kj}\hs_{ik}\frac{(\th_{ij}+1)\ttphi_{ki}}{(\th_{ik}-1)(\th_{jk}-1)}
-\sum_{a:a\not=i,j,k}\hs_{aj}\hs_{ia}\frac{\th_{ij}+1}{(\th_{ik}-1)(\th_{jk}-1)}
\cdot\frac{\ttphi_{ai}}{\th_{ka}+1}\ .\end{align}

\paragraph{4a.} The relations of the weight zero (the type 4) are also divided into 2 groups.
This is the first group of the relations:
\be \label{relation4a}[\tt_i,\tt_j]=0\ee
Note that the relations \rf{relation4a} hold for the diagonal reduction algebra
for an arbitrary reductive
Lie algebra: the images of the generators, corresponding to the Cartan sub-algebra, commute.

\paragraph{4b.}  The second group of the relations of the type 4 is (here $i\not=j$)
\be
\label{relation4}[\hs_{ij},\hs_{ji}]=\th_{ij}-\frac{1}{\th_{ij}}(\tt_i-\tt_j)^2+
\sum_{a:a\not=i,j}\left(\frac{1}{\th_{ja}+1}\hs_{ai}\hs_{ia}-\frac{1}{\th_{ia}+1}\hs_{aj}\hs_{ja}
\right)\ .\ee
The list of relations is completed.

\vskip .3cm
 Denote by $\mathfrak{R}$  the system \rf{relation1} --
\rf{relation3}, and
\rf{3.12} -- \rf{relation4} of the relations.

\begin{theorem}\hspace{-.2cm}. The relations $\mathfrak{R}$ are the defining relations for
the weight generators $\s_{ij}$ and
$t_i$ of the algebra $\Z_n$. In particular, the set  \rf{not6} of ordering relations follows over
 $\Uh$ from (and is equivalent to)  $\mathfrak{R}$.
\label{mthe}\end{theorem}

The derivation of the relations is given in \cite{KO2}; the proof of Theorem is in Section \ref{sectionproofs}.

\vskip .2cm
The relations \rf{relation1},  \rf{relation2} (a straightforward verification), as well as
\rf{relation3},  \rf{3.12},  \rf{relation4a}
and \rf{relation4}, have coefficient-bounded
terms with respect to the generators
$\hs_{ij}$ and $\tt_i$; there is no coefficient-boundedness with respect to the original
generators $\s_{ij}$ and $t_i$.
We think that the set  \rf{not6} of ordering relations is not coefficient-bounded.
\subsection{Stabilization} \label{subsection3.4}

Consider an embedding of $\gl_n$ to $\gl_{n+1}$, given by an assignment $e_{ij}\mapsto e_{ij}$,
$i,j=1,\ldots,n$, where $e_{ij}$ in the source are the generators of $\gl_n$ and target $e_{ij}$
are in $\gl_{n+1}$. The same rule $E_{ij}\mapsto E_{ij}$ defines an embedding of the Lie algebra
$\gl_n\oplus\gl_n$ to the Lie algebra $\gl_{n+1}\oplus\gl_{n+1}$ and of the enveloping algebra
$\Ar_n=\U(\gl_n\oplus\gl_n)$ to $\Ar_{n+1}=\U(\gl_{n+1}\oplus\gl_{n+1})$. This embedding clearly
maps nilpotent sub-algebras of $\gl_n$ to the corresponding nilpotent sub-algebras of $\gl_{n+1}$
and thus defines an embedding $\iota_n:\Z_n\to \Z_{n+1}$ of the corresponding double coset spaces.
However, the map $\iota_n$ is not a homomorphism of algebras. This is because the multiplication
maps are defined with the help of projectors, which are different for $\gl_n$ and $\gl_{n+1}$.

\vskip .2cm
Nevertheless,  there is an important connection between the
two multiplication maps. Namely, let $\V_{n+1}$ be
the left ideal of the algebra $\Z_{n+1}$, generated by elements $\s_{i,n+1}$, $i=1,\ldots,n$, and $\V_{n+1}'$
be the right ideal of the algebra $\Z_{n+1}$, generated by elements $\s_{n+1,i}$, $i=1,\ldots, n$. For a moment denote by
$\mmult_{(n)}:\Z_n\otimes\Z_n\to\Z_n$ and $\mmult_{(n+1)}:\Z_{n+1}\otimes\Z_{n+1}$ the multiplication maps in $\Z_n$ and
$\Z_{n+1}$ (instead of the default notation $\mult$, see \rf{not5a}).

\vskip .2cm
Let $\pi_{n+1}:\Z_{n+1}\to \Z_{n+1}$ be any linear operator in $\Z_{n+1}$, which
projects $\Z_{n+1}$ onto $\iota_n(\Z_n)$; assume that either the ideal $\V_{n+1}$
or the ideal $\V'_{n+1}$ is in the kernel of $\pi_{n+1}$,
$$\pi_{n+1}(x)=x\ , \quad x\in\iota_n(\Z_n)\ ,\qquad\text{and}\quad
\pi_{n+1}(\V_{n+1})=0\quad\text{or}\quad
\pi_{n+1}(\V'_{n+1})=0\ .$$
Define a map $\tilde{\mult}_{(n)}: \iota_n(\Z_n)\otimes \iota_n(\Z_n)\to\iota_n(\Z_n)$ as a composition
$$\tilde{\mmult}_{(n)}=\pi_{n+1}\mmult_{(n+1)}\ .$$
\begin{prop}\label{proposition1}\hspace{-.2cm}. We have a commutative diagram of maps
\be\iota_n \,\mmult_{(n)}=\tilde{\mmult}_{(n)}(\iota_n\otimes\iota_n)\ .\ee
\end{prop}

More precisely, for $i,j,k,l\leq n$ the difference
$\iota_n(z_{ij}\mmult_{(n)}z_{kl})-z_{ij}\mmult_{(n+1)}z_{kl}$ in $\Z_{n+1}$ can be written in the form
$\sum_{a=1}^n \s_{n+1,a}\mmult_{(n+1)}z_{i+k-j-l+a,n+1}\xi^{(a)}$, where $\xi^{(a)}\in \Uh$.

For the proof of Proposition, we need the following
\begin{lem}\label{lemma1}\hspace{-.2cm}. The left ideal of $\Z_n$, generated by all $\s_{in}$,
$i=1\lcd n-1$, consists of images in $\Z_n$ of  sums
$\sum_i X_i E_{in}$ with $X_i\in\Ab$, $i=1\lcd n-1$.

\vskip .2cm
The right ideal of $\Z_n$, generated by all $\s_{ni}$, $i=1\lcd n-1$, consists of images in
$\Z_n$ of  sums $\sum_i E_{ni}Y_i$ with $Y_i\in\Ab$,
$i=1\lcd n-1$.\end{lem}

{\it{Proof of Lemma}}. We follow the arguments used in the proof of the relations \rf{not3}.
 Present the projector $P$ as a sum of terms
$\xi e_{-\gamma_1}\cdots e_{-\gamma_m}e_{\gamma'_{1}}\cdots e_{\gamma'_{m'}}$,
where $\xi\in\Uh$, $\gamma_1,\dots ,\gamma_m$ and
$\gamma'_1,\dots ,\gamma'_{m'}$ are positive roots of $\f$.
{}For any $\lambda\in Q_+$ denote by $P_\lambda$  the sum of above elements
with
$\gamma_1+\dots +\gamma_{m}=
\gamma'_1+\dots +\gamma'_{m'}=\lambda$. Then $P=\sum_{\lambda\in Q_+}P_\lambda$.
{} For any $X,Y\in\Ab$ define an element  $X\mult_\lambda Y$ as the image of
 $XP_\lambda Y$ in the reduction algebra. We have
$X\mult Y=\sum_{\lambda\in Q_+} X\mult_\lambda Y$.

 For any $X\in\Ab$ and $i<n$ consider the product $X\mult_\lambda \s_{in}$. Let
 $\lambda=\sum\nolimits_{k=1}^n\lambda_k\ve_k$.
The product $X\mult_\lambda \s_{in}$ is zero if
$\lambda_n\not=0$. Indeed, in this case in each summand of $P_\lambda$ one of $e_{\gamma'_{k'}}$
is equal to some $e_{jn}$. We can order all
the monomials in $\U(\n_+)$ in such a way that all
$e_{jn}$ stand on the right. Since $[e_{jn},E_{in}]=0$, the product $e_{jn}E_{in}$ belongs
to the left ideal $\Ib$ and thus $X\mult_\lambda \s_{in}=0$
in $\Z_n$. If $\lambda_n=0$, then by PBW arguments, $P_\lambda$
 can be written
as a sum of monomials composed of generators $e_{ij}$,
$1\leq i<j<n$,  and thus their adjoint action leaves
the space, spanned by all $E_{in}$, $i<n$, invariant, so
$X\mult_\lambda \s_{in}$ is presented as an image of the sum $\sum_jX_j E_{jn}$ with $X_j\in\Ab$,
$j<n$. Thus, the left ideal, generated by $\s_{in}$ is
contained in the vector space of images in $\Z_n$ of  sums $\sum_i X_i E_{in}$.

\vskip .2cm
\noindent Moreover, $X\mult \s_{in}$ is the image of $XE_{in}+\sum_{m<i}X^{(m)}E_{mn}$ for some $X^{(m)}$
and the induction on $i$ proves the inverse inclusion.

\vskip .2cm
The second part of lemma is proved similarly.\hfill $\qed$

\vskip .2cm
{\it Proof\,} of Proposition \ref{proposition1}.
It is sufficient to prove the following
statement. Suppose $X$ and $Y$ are (non-commutative) polynomials in $E_{ij}$ with $i,j\leq n$. Then the
product of $\widetilde{\hspace{.05cm}X\hspace{.05cm}}$ and $\widetilde{\hspace{.05cm}Y\hspace{.05cm}}$
in $\Z_{n+1}$ coincides with the image in $\Z_{n+1}$ of $X\, P_{n} Y$, where $P_{n}$ is the
projector for $\gl_{n}$, modulo the left ideal in $\Z_{n+1}$, generated by all $\s_{i,n+1}$,
 $i\leq n$.
Again we note that due to the structure of the projector for any $\lambda=\sum_k\lambda_k\ve_k$ with $\lambda_{n+1}=0$,
the product $X\mult_\lambda Y$ related to $\gl_n$ coincides with product $X\mult_\lambda Y$
 related to $\gl_{n+1}$.
Thus it remains to prove that for any $X$ and $Y$ as above the element $\widetilde{\hspace{.05cm}X\hspace{.05cm}}\mult_\lambda \widetilde{\hspace{.05cm}Y
\hspace{.05cm}}$
belongs to the ideal in $\Z_{n+1}$, generated by all $\s_{i,n+1}$, $i\leq n$, once $\lambda_{n+1}\not=0$. But for
$\lambda$ with $\lambda_{n+1}\not=0$ we see, by weight arguments, that $\widetilde{\hspace{.05cm}X\hspace{.05cm}}\mult_\lambda \widetilde{\hspace{.05cm}Y
\hspace{.05cm}}$
can be presented as an image in
$\Z_{n+1}$ of the sum $\sum X_i Y_i$, such that the $(n+1)$-st component of the weight of each $Y_i$ is
not zero. Thus each $Y_i$ necessarily belongs to the left ideal generated by $E_{j,n+1}$, $j=1,\dots ,n$.
 Finally we apply Lemma \ref{lemma1} to complete the proof.

\vskip .2cm
The statement of Proposition \ref{proposition1} concerning the ideal $V_{n+1}'$ is proved similarly. \hfill $\qed$

\begin{cor}\hspace{-.2cm}.
The coefficients in the relations \rf{relation1} -- \rf{relation3},
\rf{3.12},  -- \rf{relation4} are
stable with respect to the above inclusions of $\Z_n$ to $\Z_{n+1}$.
\end{cor}
The stability of the coefficients is understood in the following sense.
Let $\cal{R}$ be a relation for $\Z_{n+1}$ from our defining list $\mathfrak{R}$,
see Subsection \ref{section3.3}.

 Assume that ${\cal{R}} $ does not contain any term with $\s_{i,n+1},\ i=1,\ldots n,$ as a left factor.
Then if we suppress in $\cal{R}$ terms which contain $\s_{i,n+1}$, $i=1,\ldots n$, as a right
factor (such term automatically contains $\s_{n+1,j}$, $j=1,\dots, n$, as a left factor), we get a relation in
$\Z_{n}$.

\vskip .2cm
Call "cut" the result of this procedure of getting the relations in $\Z_{n}$ from the relations in $\Z_{n+1}$ (under the formulated conditions).
Then all relations in $\Z_{n}$ can be obtained by cutting appropriate relations in $\Z_{n+1}$.

\vskip .2cm
Moreover, each relation in $\Z_{n}$ extends uniquely to a relation in $\Z_{n+1}$
from which it can be obtained by the cut procedure; in other words, there is a bijection between the set of relations in $\Z_{n}$ and
the set of those relations in $\Z_{n+1}$ which do not contain any term with $\s_{i,n+1}$, $i=1,\ldots n$, as a left factor.

\vskip .2cm
The stabilization rule is certainly not an isolated $\gl$ phenomenon; it can be generalized to certain other quadruplets of algebras
replacing those which participate in the diagram
$${{\gl_n}^{\textstyle{\nearrow}}_{\textstyle{\searrow}}}\begin{array}{c}\gl_n\oplus\gl_n\\[1em]\gl_{n+1}\end{array}
\! \phantom{}^{\textstyle{\searrow}}_{\textstyle{\nearrow}} \gl_{n+1}\oplus\gl_{n+1}\ .$$

\setcounter{equation}0
\section{Completeness of relations\vspace{.25cm}}\label{sectionproofs}

\paragraph{1.} We first give general arguments, proving the weakened version, in which $\Uh$ is
enlarged to $\UUh$, of Theorem \ref{mthe}.

\vskip .2cm As before, denote by $\mathfrak{R}$  the system
\rf{relation1},  \rf{relation2}, \rf{relation3}, \rf{3.12},
\rf{relation4a} and \rf{relation4} of relations. We shall see that
it is equivalent to the system  \rf{not6} of the ordering rules. The
system $\mathfrak{R}$ follows from \rf{not6} since \rf{not6} is the
set of defining relations for the weight generators; we have to
verify the opposite implication. For a moment denote the generators
from the set $\{\hs_{ij},\tt_i\}$ by symbols $\pp_L$, labeled by a
single index $L$, $L=1,2,\dots ,n^2$. The number of ordering rules
for $n^2$ variables $\pp_L$ is $n^2(n^2-1)/2$. So, to prove the
completeness, it is sufficient to show that the dimension of the
subspace (over  $\UUh$) spanned by $\mathfrak{R}$ is at least
$n^2(n^2-1)/2$. Any relation from $\mathfrak{R}$ is a sum of
products $\pp_L\mult\pp_M$ with coefficients in $\Uh$ plus,
possibly, a term of zero degree in $\pp$'s. Denote by
$\mathfrak{R}_0$  the system $\mathfrak{R}$  with degree zero terms
dropped. It suffices to show that \be \mathrm{the}\
\mathrm{system}\ \,\mathfrak{R}_0\,\ \mathrm{contains}\,\
{n^2(n^2-1)}/{2}\ \, \mathrm{linearly}\
\mathrm{independent}\ \mathrm{over}\,\ \UUh\ \mathrm{relations}\ .\label{fafa}\ee Once the
coefficients from $\Uh$ in all relations from $\mathfrak{R}_0$ are
placed on the same side, say, on the right, from the monomials
$\pp_L\mult\pp_M$, one can give arbitrary numerical values to the
variables $\th_{ij}$  (respecting  linear dependencies between them).
 To check the assertion \rf{fafa} it is enough
to find a set of values for which the corresponding system with
numerical coefficients has $n^2(n^2-1)/2$ linearly independent
relations . But when all $\th_{ij}$ tend to $\infty$ (in the
following way: $\th_{i,i+1}=c_{i,i+1}h$, $h\rightarrow\infty$ and
$c_{i,i+1}$ are constants),
 we directly observe that the system
$\mathfrak{R}_0$ becomes simply $\pp_L\mult\pp_M=\pp_M\mult\pp_L$, $M>L$. The proof of the
completeness over $\UUh$ is finished. \hfill $\qed$

Note that we did not use in the above arguments the compatibility of the ordering $\prec$ with the partial
order $<$ on $\h^*$.

\paragraph{2.} Given an order, let $X$ be a formal vector of all unordered products $\pp_L\mult\pp_K$ and $Y$ a formal vector of all ordered products.
To rewrite $\mathfrak{R}$ in the form of ordering relations, one has to solve for $X$ a linear system of equations
\be {\cal{A}}X={\cal{B}}Y+C\ ,\label{caode}\ee
where $C$ is a vector of degree 0 terms; ${\cal{A}}$ and ${\cal{B}}$ are certain matrices with coefficients in $\Uh$ (by the above proof,
${\cal{A}}$ is a square matrix). The solution of this system may cause an appearance of coefficients from $\UUh$ (not from $\Uh$) in the ordering
relations. This happens, for example, for the lexicographical order for the generators $z_{ij}$ and $t_i$ (with $z_{ii}=t_i$) of $\Z_n$ for $n>2$
(we don't give details; it is an explicit calculation).
 It follows from  \rf{not3} (and statement (e) of Section \ref{section2new}) that
 for the order \rf{not4} - \rf{not4b} the solution of the system \rf{caode} is defined over
 the ring $\Uh$. However, this shows only that for this order possible terms from $\UUh$ in the determinant of ${\cal{A}}$ simplify in the combinations
${\cal{A}}^{-1}{\cal{B}}$ and ${\cal{A}}^{-1}C$; the systems $\mathfrak{R}$ and
$\mathfrak{R}^\prec$ may still be not equivalent over $\Uh$, in the sense that the
elements of the matrix ${\cal{A}}^{-1}$ may not belong to $\Uh$
  and we cannot transform the system $\mathfrak R$ to the system  $\mathfrak{R}^\prec$ by composing
  linear over $\Uh$ combinations of relations from $\mathfrak R$.

\paragraph{3.} We now pass to the proof of Theorem \ref{mthe}.
Let ${\cal{F}}$ be the free algebra with the weight generators $\s_{ij}$ and $t_i$ over $\Uh$. Let $\mathfrak{R}^\prec$ be the set of
ordering relations \rf{not6}. Both $\mathfrak{R}$ and $\mathfrak{R}^\prec$ are defined over $\Uh$ and we have the homomorphism
$\varpi :\, {\cal{F}}/\mathfrak{R}\rightarrow {\cal{F}}/\mathfrak{R}^\prec$. According to the weak form of theorem \ref{mthe}, see paragraph 1 of this
subsection, the homomorphism $\varpi$ becomes the isomorphism after taking the tensor product with $\UUh$ (over $\Uh$). We shall now prove
that $\varpi$ itself is the isomorphism.

 The proof is done by induction in $n$, with the help of the stabilization law and an explicit
  calculation of certain determinants
 (and one can follow the precise structure of appearing denominators at each step).
The induction base is $n=1$, there is nothing to prove for $\Z_1$.

\vskip .2cm
All we have to show in general case is that the numerator of the determinant of the matrix
${\cal{A}}$, figuring in \rf{caode}, is a product of linear
factors of the form
\rf{musel}. The relations are weighted so the matrix ${\cal{A}}$ has a block structure, blocks ${\cal{A}}_\lambda$ are labeled by the relation weights.
The determinant of ${\cal{A}}$ is the product of the determinants of the blocks ${\cal{A}}_\lambda$.

  Consider $\Z_{n-1}$ as a subspace in $\Z_n$ as in section \ref{subsection3.4}. Fix a weight $\lambda$ for $\Z_{n-1}$. Call
${\cal{L}}_\lambda^{(n)}$ the linear subsystem ${\cal{A}}_\lambda X_\lambda ={\cal{B}}_\lambda Y_\lambda +C_\lambda$ of \rf{caode}, corresponding
to the weight $\lambda$ for $\Z_n$. The system ${\cal{L}}_\lambda^{(n)}$ contains the subsystem $_{\phantom{\lambda}}^{(n)}
\! {\cal{L}}_\lambda^{(n-1)}$, corresponding to the generators from $\Z_{n-1}$ (recall that the relations are labeled by pairs of generators, so the
subsystem  $_{\phantom{\lambda}}^{(n)}\! {\cal{L}}_\lambda^{(n-1)}$ is well defined). Compare $_{\phantom{\lambda}}^{(n)}\! {\cal{L}}_\lambda^{(n-1)}$
with the corresponding system ${\cal{L}}_\lambda^{(n-1)}$ for $\Z_{n-1}$. By the stabilization principle, the system ${\cal{L}}_\lambda^{(n-1)}$ is the
cut of the system $_{\phantom{\lambda}}^{(n)}\! {\cal{L}}_\lambda^{(n-1)}$ in the sense of section \ref{subsection3.4}: there is a bijection
between the
two systems and the relations from $_{\phantom{\lambda}}^{(n)}\! {\cal{L}}_\lambda^{(n-1)}$ have, compared to the corresponding relations from
${\cal{L}}_\lambda^{(n-1)}$, extra terms with $\s_{ni}\mult \s_{jn}$ for certain $i,j<n$. By induction, ${\cal{L}}_\lambda^{(n-1)}$ is equivalent, over
its own $\Uh$, to the system of ordering relations. Making the same transformation with the system $_{\phantom{\lambda}}^{(n)}
\! {\cal{L}}_\lambda^{(n-1)}$ preserves the ordered form since the terms $\s_{ni}\mult \s_{jn}$ are ordered. This argument shows that we need to
consider only the subset of relations labeled by those pairs of generators $(\pp_L,\pp_M)$ for which $\pp_L$ or $\pp_M$ do not belong to $\Z_{n-1}$,
\be (\pp_L,\pp_M)\, :\ \pp_L\notin \Z_{n-1}\ {\mathrm{or}}\ \pp_M\notin \Z_{n-1}\ .\label{indrests}\ee
Applying the just constructed ordering rules (equivalent to the system $_{\phantom{\lambda}}^{(n)}\! {\cal{L}}_\lambda^{(n-1)}$) to these
remaining relations, we leave in them only ordered terms $\pp_{L'}\pp_{M'}$, $L'<M'$,
 with two generators from $\Z_{n-1}$, $\pp_{L'},\pp_{M'}\in \Z_{n-1}$.
\vskip .2cm

  We shall now consider separately each type of weight relations  listed
in the beginning of  Section \ref{section3.3}.
The relations of types 1 and 2 do not cause any difficulty.

\vskip .2cm
The number of relations of the types 3 or 4 grows with $n$.  The change of variables \rf{3.1a}, as well as the renormalization \rf{not8}
and its inverse,
have allowed denominators, so we can work with the generators $\tt_i$ and $\hs_{ij}$ instead of $t_i$ and $\z_{ij}$.

\paragraph{4. Relations  of type 4.} For  the relations \rf{relation4a} and  \rf{relation4}, the restriction \rf{indrests} shows that we have to
consider only the subsystem, corresponding to the pairs $(\z_{ni},\z_{in})$, $i<n$, of the relations  \rf{relation4}.
By the arguments from the paragraph
above, we assume that the only unordered quadratic monomials in this subsystem are
 $\hs_{in}\mult\hs_{ni}$, $i<n$. Rewrite this subsystem in the form
\rf{caode}:
\be \hs_{in}\mult\hs_{ni}+\sum_{a:a<n,a\neq i}\frac{1}{\th_{ia}+1} \hs_{an}\mult\hs_{na}=\dots\ ,\ee
where dots stand for ordered terms. Therefore, the matrix $\mathring{\mathfrak{A}}$, whose determinant we need to calculate, is simply
\be \mathring{\mathfrak{A}}_{ij}:=\frac{1}{\th_{ij}+1}\ ,
\label{bacd}\ee
where, we recall, $\th_{ij}=\th_{i}-\th_{j}$; in particular, $\th_{ii}=0$. The determinant of such matrix is well known.
The matrix $\mathring{\mathfrak{A}}$ is
the specialization of the matrix
\be {\mathfrak{A}}_{ij}:=\frac{1}{x_i+y_j}
\label{cma}\ee
at $x_i=\th_i$ and $y_j=-\th_j+1$. The determinant of the matrix $ {\mathfrak{A}}$, calculated in \cite{Ca}, is $\det {\mathfrak{A}}=\prod_{i,j:i<j}\,
\Bigl( (x_i-x_j)(y_i-y_j)\Bigr)/\prod_{i,j}\, (x_i+y_j)$. It follows that $\det \mathring{\mathfrak{A}}=\prod_{i,j:i<j}\, \th_{ij}^2/(\th_{ij}^2-1)$.
The inverse of
$\mathring{\mathfrak{A}}$ has thus allowed denominators.

\paragraph{5. Relations  of type 3.}   For the relations  \rf{relation3} and  \rf{3.12} of the type 3, the restriction \rf{indrests}
shows that we have to
consider only the relations of the weights $\ve_i-\ve_n$ and $\ve_n-\ve_i$, $i<n$.

\vskip .2cm
We start with the weight $\ve_i-\ve_n$ with a fixed $i$, $i<n$. The unordered quadratic monomials of the weight $\ve_i-\ve_n$ are
 \begin{eqnarray}&\hs_{an}\mult\hs_{ia}&\ \ \ \ {\mathrm{with}}\ \ a:\, 2a<i+n\ ,\ a\neq i\ ,\label{fcf31}\\[.5em]
 &\hs_{ij}\mult\hs_{jn} &\ \ \ \ {\mathrm{with}}\ \ j:\, i+n\leq 2j\ ,\label{fcf32}\\[.5em]
&\hs_{in}\mult\tt_b\ .&\label{fcf33}\end{eqnarray}
All relations \rf{relation3} and  \rf{3.12} participate in our system. However, the system is block-triangular and can be analyzed.

\vskip .2cm
Denote by ${\mathfrak{r}}_{ikl}$ the relation from the list \rf{relation3} whose left hand side starts with $\hs_{ik}\mult\hs_{kl}$.  Let
$\boldsymbol{\kappa}_{ikl}$
be the
coefficient of the term $\hs_{ik}\mult\hs_{kl}$ in ${\mathfrak{r}}_{ikl}$. The relations ${\mathfrak{r}}_{ijn}$ can be rewritten in the form
(since $\th_{jj}=0$)
\be \hs_{ij}\mult\hs_{jn}\,\boldsymbol{\kappa}_{ijn}=\sum_{a:a\neq i,a<n}\hs_{an}\mult\hs_{ia}\frac{B_{ai}}{\th_{ja}+1}+\dots\ .\ee
Here dots stand for ordered terms with $\tt_b\mult\hs_{in}$
(the term with $\hs_{jn}\mult \hs_{ij}$ is absorbed into the sum). Among the unordered
monomials \rf{fcf31}-\rf{fcf33} only the monomials \rf{fcf31} enter the relations ${\mathfrak{r}}_{ijn}$ with $j$ such that $2j<i+n$ and $j\neq i$.
Thus the
subsystem $\{ {\mathfrak{r}}_{ijn}\, |\, j:\, 2j<i+n\, ,\, j\neq i\}$ contains as many relations as unordered monomials.
The matrix, whose determinant we have
to calculate in order to express,  using this subsystem, the unordered monomials \rf{fcf31} in terms of ordered monomials is
$\mathring{\mathfrak{A}}_{aj}':=
B_{ai}/(\th_{ja}+1)$; the $a$-th row contains $B_{ai}$ as  the common factor, so the determinant of the matrix $\mathring{\mathfrak{A}}'$ is the
product of $B_{ai}$ (over $a$ such that $\, 2a<i+n$ and $a\neq i$) times the determinant
 of the matrix of the same form \rf{bacd} as before. Thus the inverse of the
determinant of   the matrix $\mathring{\mathfrak{A}}'$  belongs to $\Uh$.
 We use this subsystem to order the monomials \rf{fcf31}.

\vskip .2cm
After the monomials \rf{fcf31} are ordered, the rest of the relations ${\mathfrak{r}}_{ijn}$ (with $j:\, i+n\leq 2j$) turns into the set of
the ordering
relations for the monomials \rf{fcf32}; each relation contains exactly one unordered monomial of the form  \rf{fcf32} with the coefficient
 $\boldsymbol{\kappa}_{ijn}$ whose inverse has allowed determinants.

\vskip .2cm
The set of relations \rf{3.12} provides the ordering rules for the monomials $\hs_{in}\mult\tt_k$ once one knows the ordered expressions
for all monomials
$\hs_{an}\mult\hs_{ia}$.

\paragraph{6. Relations  of type 3, weight $\ve_n-\ve_i$.} The considerations of
Section \ref{section2new} show that for any two orders on the weight basis of $\p$,
 compatible with the partial order $<$ on $\h^*$, the ordering relations \rf{not3}
 for them  are
equivalent over $\Uh$. Define, instead of \rf{not4}- \rf{not4b}, the order $\grave{\prec}$ by
\begin{equation}\label{not4p}\z_{ij}\grave{\prec}\z_{kl}\quad \text{if}
\quad i-j>k-l\quad\text{or}\quad \left\{\begin{array}{ll}i>k\quad &{\mathrm{if}}
\quad i-j=k-l>0\ ,\\[.3em]
i<k\quad &{\mathrm{if}}\quad i-j=k-l<0\ ,\\[.3em] {\mathrm{arbitrarily}}
\quad &{\mathrm{if}}\quad i-j=k-l=0\ .\end{array}\right.\end{equation}
The peculiarity of the order $\grave{\prec}$ is that the anti-involution
$\epsilon$, see \rf{anep}, transforms the set of quadratic ordered monomials
of any non-zero weight $\lambda$, $\lambda\neq 0$, into the set of quadratic
ordered monomials of the weight $(-\lambda)$.

It is proved in \cite{KO2}, that the system $\mathfrak{R}$ is closed under the
anti-involution $\epsilon$  (that is,  $\mathfrak{R}$ and
$\epsilon (\mathfrak{R})$ are equivalent over $\Uh$).
{}For the order $\grave{\prec}$, the application of the anti-involution
$\epsilon$ reduces the question about the equivalence over
$\Uh$ of $\mathfrak{R}$
and the set \rf{not6} of the ordering relations for the weight $\ve_n-\ve_i$
to the same question
for the weight $\ve_i-\ve_n$.
By the preceding paragraph,
the equivalence assertion follows for the order $\grave{\prec}$ and therefore for any other order,
compatible with the partial order $<$ on $\h^*$, for
example, the order $\prec$.

The proof of the theorem \ref{mthe} is completed. \hfill $\qed$

\vskip .2cm
The set $\mathfrak{R}^\prec$ of ordering relations \rf{not6} is, by construction, closed over
$\Uh$ under the involution $\omega$, see \rf{not2a}. As a
by-product of the equivalence of $\mathfrak{R}$ and $\mathfrak{R}^\prec$ over $\Uh$ we observe
that $\mathfrak{R}$ is closed over $\Uh$ under the
involution $\omega$ as well.

\vskip .2cm
Note that all denominators, which appeared in the proof, are of the form $\th_{ij}\pm \varsigma$,
$i<j$, where $\varsigma =0,1$ or 2.

\paragraph{7.} As the proof shows, essentially the only  matrix we have to invert is of the
 form \rf{bacd}. The matrix, inverse to \rf{bacd} reads
\be  (\mathring{\mathfrak{A}}^{-1})_{ij}=-\frac{1}{\th_{ij}-1}\prod_{a:a\neq i}
\frac{\th_{ia}-1}{\th_{ia}}
\prod_{b:b\neq j}\frac{\th_{jb}+1}{\th_{jb}}\ .\label{bacdinv}\ee
The verification of \rf{bacdinv} in the form $\sum_j (\mathring{\mathfrak{A}}^{-1})_{ij}
 \mathring{\mathfrak{A}}_{jk}=\delta_{jk}$,
where $\delta_{jk}$ is the
Kronecker delta, reduces to the identity
\be\frac{1}{\th_{ik}+1}\prod_{b:b\neq i} \frac{\th_{ib}+1}{\th_{ib}}-\sum_{j:j\neq i}
\frac{1}{\th_{ij}(\th_{jk}+1)}
\prod_{b:b\neq i,j} \frac{\th_{jb}+1}{\th_{jb}}=\delta_{ik}\prod_{b:b\neq i}
\frac{\th_{ib}}{\th_{ib}-1}\ ,\ee
which is checked by an evaluation of residues and the values at infinity of
both sides as functions of $\th_i$.

\vskip .2cm
The inverse of the more general matrix \rf{cma} reads
\be  ({\mathfrak{A}}^{-1})_{ij}=(x_j+y_i)\prod_{a:a\neq j}\frac{x_a+y_i}{x_a-x_j}
\prod_{b:b\neq i}\frac{y_b+x_j}{y_b-y_i}\ .\label{cmainv}\ee
It is demonstrated similarly to \rf{bacdinv}, by an appropriate evaluation of residues
and the values at infinity.

\vskip .2cm
The formula \rf{bacdinv} is equivalent (not directly equal) to the specialization of \rf{cmainv}
at $x_i=\th_i$ and $y_j=-\th_j+1$.

\vskip .2cm
The formula \rf{bacdinv} provides a recursive way to transform the system $\mathfrak{R}$
into the set of ordering relations.

\section*{Acknowledgments}

We thank Lo\"{\i}c Poulain d'Andecy for an independent partial confirmation of correctness of the relations from subsection
\ref{section3.3} by a computer aided check of the Poincar\'e--Birkhoff--Witt theorem for the  algebra $\Z_4$ defined by the
generators and relations. We are indebted to Elena Ogievetskaya for valuable help in preparation of the manu\-script.

\vskip .2cm A part of the present work was done during visits of S.
K. to CPT and CIRM in Marseille. The authors thank the staff of the
Institutes for providing excellent working condi\-tions during these
visits. S. K. was supported by the RFBR grant 08-01-00667, joint
 CNRS-RFBR grant 09-01-93106 and
the grant for Support of Scientific Schools 3036-2008-2. Both
authors were supported by the ANR project GIMP No.
ANR-05-BLAN-0029-01.

\end{document}